\documentclass{article}
\usepackage{xcolor}
\usepackage{enumitem}
\usepackage{comment}
\usepackage{mathtools, nccmath}
\usepackage{amsmath,amssymb,amsthm}
\usepackage{esvect}
\usepackage{graphicx}
\usepackage{hyperref}
\hypersetup{
    colorlinks=true,
    citecolor=blue,
    linkcolor=blue,
    filecolor=blue,      
    urlcolor=blue,
}
\usepackage{float}
\graphicspath{ {./images/} }
\usepackage{placeins}

\providecommand{\keywords}[1]
{
  \small	
  \textbf{\textit{Keywords---}} #1
}
\newtheorem{definition}{Definition}[section]
\newtheorem{lemma}[definition]{Lemma}
\newtheorem{remark}[definition]{Remark}
\newtheorem{proposition}[definition]{Proposition}
\newtheorem{theorem}[definition]{Theorem}
\newtheorem{corollary}[definition]{Corollary}
\newtheorem{example}[definition]{Example}

\newcommand{\dom}{\operatorname{dom}}
\newcommand{\R}{{\rm I\!R}}
\newcommand{\ones}{{\mathbf{1}}}

\newcommand{\N}{\mathbb{N}}

\newcommand{\od}{\overline{\Delta}}
\newcommand{\de}{\Delta}
\newcommand{\bbm}{\begin{bmatrix}}
\newcommand{\ebm}{\end{bmatrix}}
\newcommand{\ds}{\delta_f}
\newcommand{\limn}{\lim_{N \to \infty}}
\newcommand{\limnvec}{\lim_{\vv{N} \to \vv{\infty}}}

\newcommand{\nintrecc}{\int_{R(0;d)}}

\newcommand{\dexxe}{\, dx}

\newcommand{\intballzero}{\int_{B_n(0;r)}}

\DeclareMathOperator{\id}{Id}

\DeclareMathOperator{\diag}{diag}

\begin{document}

\title{Limiting behaviour of the generalized simplex gradient as the number of points tends to infinity on a fixed shape in $\R^n$}

\author{
{\sc Warren Hare\thanks{Corresponding author. Email: Warren.Hare@ubc.ca}} \\[2pt]
Department of Mathematics,\\ University of British Columbia, Okanagan campus,\\ 3187 University Way, Kelowna, BC, Canada,\\
{\sc and}\\[6pt]
{\sc Gabriel Jarry-Bolduc}\thanks{Email: gabjarry@alumni.ubc.ca}\\[2pt]
Department of Mathematics, \\University of British Columbia, Okanagan campus \\
{\sc and}\\[6pt]
{\sc Chayne Planiden}\thanks{Email: chayne@uow.edu.au }\\[2pt]
School of Mathematics and Applied Statistics,\\ University of Wollongong, Wollongong, NSW, 2500, Australia.
}

\maketitle

\begin{abstract}
This work investigates the asymptotic behaviour of the gradient approximation method called the generalized simplex gradient (GSG). This method has an error bound that at first glance seems to tend to infinity as the number of sample points increases, but with some careful construction, we show that this is not the case. For functions in finite dimensions, we present two new error bounds ad infinitum depending on the position of the reference point. The error bounds  are not a function of  the number of sample points and thus remain finite.
\end{abstract}


\section{Introduction}
In optimization, it is common to find oneself in a setting in which it is necessary or desirable to approximate gradients of functions \cite{BBN2018,conn2008geometry,Maggiar2018,MMSMW2017,powell2009bobyqa,verderio2017construction}. This occurs most commonly in Derivative-Free Optimization (DFO), where gradients are computationally expensive to calculate, inconvenient to use or simply impossible to obtain \cite{AudetHare_DFObook,conn-scheinberg-vicente-2009}. 

One of the main concerns regarding any DFO gradient approximation method is the establishment of an upper bound on the error between the gradient approximation and the true gradient of the function. If the error can be controlled, then efficient, convergent optimization algorithms can be designed around the method. Recent work \cite{hare2020error} has been done in this regard on several expansions of a popular approximation technique called the simplex gradient. This paper focuses on the expansion known as the \emph{generalized simplex gradient}. The definition of this method is given in the next section.

Given a $\mathcal{C}^{2}$ function $f:\R^n\to\R^p$, the simplex gradient approximates the gradient $\nabla f$ at a point $x$ by building an affine function $\varphi$ that is a good approximation to $f$ near $x$ and then calculating $\nabla\varphi$. The construction is accomplished by selecting $n$ suitably spaced points in the domain of $f$ that together with the reference point $x$ form a simplex and serve to define $\varphi$. Then $\nabla\varphi$ is calculated using $n+1$ function evaluations of $f$. Some early examples of the simplex gradient are found in \cite{MR323040,MR0326978}, in which a particular simplex based on coordinate directions is used to approximate gradients of quadratic and cubic functions. The {\em Simplex Gradient} is introduced in \cite{Kelley99Book}, where the ideas are generalized to allow for an arbitrary simplex.  The error bound on the simplex gradient is a function of the dimension $n$ and the radius of the simplex $\Delta_S$ \cite{Kelley99Book}. 

The {\em Generalized Simplex Gradient} (GSG) relaxes the requirement of evaluating exactly $n+1$ points in $\R^n$, making it possible to obtain an approximate gradient with controlled error bound by using any finite number of points $m+1$ (the positive integer $m$ can be greater than, equal to, or less than $n$) \cite{Regis2015}. This method has an error bound that depends on the number of points used $m$ and the sampling radius $\Delta_S$, which means that if one were to increase the number of sample points used and consider the limit as $m\to\infty$, then the error bound does not necessarily remain bounded.  Indeed, it is possible to construct examples in which the classical error bound tends to infinity as $m$ increases \cite{MR4163088}. This is a counter-intuitive event, as it is reasonable to conjecture that more sample points would provide better accuracy of the model function. This problem is investigated in \cite{MR4163088} for functions $f:\R\to\R$, in which new error bounds are developed and shown to have a more desirable behaviour at the limit.

In  \cite{MR4163088}, it is shown that the norm of the difference between the approximate derivative and the true derivative of a one-dimensional function has a limit that is a factor of the Lipschitz constant of $\nabla f$ and the sampling radius $\Delta_S$. However, \cite{MR4163088} considers only single-variable functions. The purpose of the present work is to extend those results to the multivariable setting. 

We first consider sampling over a hyperrectangle with the reference point at a corner point.  We explore the limiting behaviour of the gradient approximation on $\R^n$ and calculate the limit of the corresponding error bound as the number of sample points tends to infinity.  We then repeat the process considering the case where the sampling set is a ball with the reference point at the center.   From these results, it becomes clear how the method can be adapted to other shapes. We discuss the option of other shapes further in the concluding section.

The remainder of this paper is organized as follows. Section \ref{sec:Prel} defines the notation used throughout the paper and provides some definitions for later use. In Section \ref{sec:hyperrectangle}, we investigate the limiting behaviour of the GSG in the Cartesian setting when the sample region is a hyperrectangle, and in Section \ref{sec:error}, we present its error bound ad infinitum. Section \ref{sec:ball} considers a ball as the sample region and provides an error bound ad infinitum.  Illustrative examples are included throughout.   We make some concluding remarks and discuss avenues of future research in Section \ref{sec:conc}.

\section{Preliminaries} \label{sec:Prel}

\subsection{Notation}

Throughout this work, we use the standard notation found in \cite{rockwets}. The domain of a function $f$ is denoted by $\dom f$. The transpose of a matrix $A$ is denoted by $A^\top$. We work in finite-dimensional space $\R^n$ with inner product $x^\top y=\sum_{i=1}^nx_iy_i$ and induced norm $\|x\|=\sqrt{x^\top x}$. The identity matrix in $\R^{n \times n}$ is denoted by $\id_n$, or by $\id$ if the dimension is clear. The vector $e_i \in \R^n$ denotes the $i$\textsuperscript{th} column of $\id_n.$ The vector of all ones in $\R^n$ is denoted by $\ones_n$, or $\ones$ if the dimension is clear.  The $i$\textsuperscript{th} component of a vector $v$ is denoted by $v_i$, and of a vector $v_n$ by $[v_n]_i$. Similarly,  the entry in the $i$\textsuperscript{th} row and $j$\textsuperscript{th} column of a matrix $A$ is denoted by $A_{i,j},$ and of a matrix $A_n$ by $[A_n]_{i,j}$. The diagonal matrix $D$ in $\R^{n \times n}$  is written $\diag\bbm D_{1,1}&D_{2,2}&\cdots&D_{n,n} \ebm$ when convenient. Given a matrix $A \in \R^{n \times m},$ we use the induced matrix norm 
\begin{align*}
    \Vert A \Vert=\Vert A \Vert_2=\max \{ \Vert Ax\Vert_2 \, : \, \Vert x \Vert_2=1 \}.
\end{align*}
The sphere and the  ball of radius $r>0$ centered at  $x^0 \in \R^n$  are denoted by $S_n(x^0;r),$ and  $B_n(x^0;r)$  respectively. That is,
\begin{align*}
    S_n(x^0;r)&=\left \{  x \in \R^n: \Vert x-x^0 \Vert=r \right \}, \quad
    B_n(x^0;r)=\left \{ x \in \R^n: \Vert x-x^0 \Vert \leq r \right \}.
\end{align*}
\subsection{Definitions and minor results}

In this section, we list some concepts and properties that will be useful in developing the main results, as well as the formal definition of the GSG. Central to the GSG is the Moore--Penrose matrix pseudoinverse.

\begin{definition}[Moore--Penrose pseudoinverse] \label{def:mpinverse}
Let $A\in\R^{n\times m}$. The \emph{Moore--Penrose pseudoinverse} of $A$, denoted by $A^\dagger$, is the unique matrix in $\R^{m \times n}$ that satisfies the following four equations:
\begin{align*}
AA^\dagger A&=A,\\A^\dagger AA^\dagger&=A^\dagger,\\(AA^\dagger)^\top&=AA^\dagger,\\(A^\dagger A)^\top&=A^\dagger A.
\end{align*}
\end{definition}

Note that for any $A \in \R^{n \times m},$ there exists a unique Moore-Penrose pseudoinverse $A^\dagger\in\R^{m\times n}.$ The following property will be used frequently in the sequel.
\begin{itemize}
 \item If $A$ has full row rank $n$, then $A^\dagger$ is a right-inverse of $A$, so that $AA^\dagger=\id_n.$ In this case, $A^\dagger=A^\top (A A^\top)^{-1}.$ 
\end{itemize}
\noindent Next, we introduce the definition of the GSG and provide its error bound.

\begin{definition}[Generalized simplex gradient] Let $f:\dom f \subseteq \R^n\to\R$.  Let $x^0 \in \dom f$ be the reference point.  Let $S \in \R^{n \times N}$ with $x^0+Se_j \in \dom f$ for all $j\in\{1,2,\ldots,N\}$. The \emph{generalized simplex gradient} of $f$ at $x^0$ over  $S$ is denoted by $\nabla_s f(x^0;S)$ and defined by
\begin{equation} \label{eq:GSG} \nabla_s f(x^0;S)=(S^\dagger)^\top \ds (x^0;S),\end{equation}
where $\ds (x^0;S)=\bbm f(x^0+Se_1)-f(x^0)& \cdots & f(x^0+Se_N)-f(x^0) \ebm^\top \in \R^N.$
\end{definition}
Occasionally, \eqref{eq:GSG} is written in terms of $(S^\top)^\dagger$.  These two forms are equivalent, as $(S^\top)^\dagger=(S^\dagger)^\top.$ The following theorem establishes an error bound for the GSG. The error bound depends on the radius $\Delta_S$ of the matrix $S$, that is,
\begin{align}
    \Delta_S&=\max \left\{ \Vert Se_j \Vert : j \in \{1, 2, \dots, N\} \right\}.\label{df:radius}
\end{align} 
\begin{theorem}[Classical error bound for the GSG]\emph{\cite[Cor.1]{Regis2015} \& \cite[Thm.3.3]{hare2020error}}\label{prop:GSGerror}
 Let $S \in \R^{n \times N}$ have full row rank and radius $\Delta_S$. Let $f$ be $\mathcal{C}^{2}$ on an open domain containing $B_n(x^0;\Delta_S)$ where $x^0$ is the reference point and $\Delta_S>0$. 
 Then
\begin{align}\label{eq:ebgsg}
\|\nabla_s  f(x^0;S)-\nabla f(x^0)\| &\leq \frac{\sqrt{N}}{2}L_{\nabla f}  \left \Vert (\widehat{S}^\top)^\dagger \right \Vert \Delta_S,
\end{align}
where $\widehat{S}=S/\Delta_S$ and $L_{\nabla f}$ denotes the Lipschitz constant of $\nabla f$ on $B_n(x^0;\Delta_S).$

Moreover, if $f \in \mathcal{C}^3$ on an open domain containing $B_n(x^0;\Delta_S)$ and $S$ (or some permutation of $S$) has the form $S=\bbm A&-A\ebm$ for some  $A \in \R^{n \times \frac{N}{2}},$  then 
\begin{align}\label{eq:ebgcsg}
\|\nabla_s  f(x^0;S)-\nabla f(x^0)\| &\leq \frac{\sqrt{N}}{6}L_{H}  \left \Vert (\widehat{A}^\top)^\dagger \right \Vert \Delta_S^2,
\end{align}
where $\widehat{A}=A/\Delta_S$ and $L_H$ denotes the Lipschitz constant of $\nabla^2 f$ on  $B(x^0;\Delta_S).$
\end{theorem}
Note that \eqref{eq:ebgcsg} is the error bound defined for the \emph{generalized centered simplex gradient (GCSG)} in \cite[Thm.3.3]{hare2020error}. It was shown in \cite[Prop.2.10]{hare2020error} that the GSG over $S=\bbm A&-A \ebm$ and the GCSG over $A$ are equivalent. So, we consider the GCSG a specific case of the GSG. 

Finally, we recall the Sherman--Morrison--Woodbury formula for calculating matrix inverses.
\begin{definition}[Sherman--Morrison--Woodbury formula]
Let $A \in \R^{n \times n}$ be a non-singular matrix and $u,v \in \R^n.$ If $1+v^\top A^{-1} u \neq 0,$ then 
\begin{equation}\label{eq:sherman}
    (A+uv^\top)^{-1}=A^{-1}-\frac{A^{-1}uv^\top A^{-1}}{1+v^\top A^{-1}u}.
\end{equation}
\end{definition}

\section{The GSG over a dense hyperrectangle}\label{sec:hyperrectangle}

In this section, we find an expression for the GSG on a hyperrectangle $R$ as the number of points in $R$ tends to infinity in such a way that they form a dense grid.

As the case $n=1$ is covered in \cite{MR4163088}, for the remainder of this paper we assume $n \geq 2.$

\subsection{Preliminaries: Using the rightmost endpoints in each partition} 

Let $x^0 \in \R^n$  be the reference point.  Consider a  hyperrectangular sample region  with side lengths $\de_1, \de_2, \dots, \de_n >0$ where $x^0$ is the leftmost vertex of the hyperrectangle (the point with the lowest value for each component of all points in the hyperrectangle). We denote this hyperrectangle by $R(x^0;d)$ where $d=\bbm \de_1 & \de_2 & \cdots & \de_n \ebm^\top.$   Then $R(x^0;d)$ is partitioned into sub-hyperrectangles with lengths $\od_i=\de_i/N_i$ where $N_i \in\N\setminus\{1\}$ for all $i \in \{1, 2, \dots, n\}.$  We define $N$ to be the product of all $N_i$:  $N=N_1N_2 \cdots N_n$. Then $S_R$ is defined to be a matrix in $\R^{n \times N}$ that  contains all the directions used to form the sample points when the rightmost  endpoint of each sub-hyperrectangle is chosen.\footnote{The process will be generalized later in this section to allow choosing any arbitrary point in each partition of $R(x^0;d)$.} Hence, $R(x^0;d)$ contains $N$ sample points and one reference point, $x^0$.  

Let $\overline{D}=\diag \bbm \od_1&\od_2&\cdots&\od_n \ebm \in \R^{n \times n}$. Then $S_R$ can be  written as a block matrix in the following way:
\begin{align}\label{eq:SR}
    S_R&=\bbm B_R^{1, 1, \dots,1,1}&B_R^{1, 1, \dots,1, 2}&\cdots&B_R^{N_2, N_3, \dots, N_{n-1}, N_n} \ebm \in \R^{n \times N}, 
\end{align}
where 
\begin{align}  \label{eq:BRabm}
B_R^{\vv{z}}&=B_R^{z_2,z_3, \dots,z_{n-1}, z_{n}}\\
&= \overline{D} \ddot{B}_R^{\vv{z}} \notag \\
&=\overline{D} \bbm 1&2&\cdots&N_1\\z_2&z_2&\cdots& z_2\\z_3&z_3&\cdots&z_3\\ \vdots&\vdots&\cdots&\vdots\\z_{n-1}&z_{n-1}&\cdots&z_{n-1}\\z_{n}&z_{n}&\cdots&z_{n} \ebm \in \R^{n \times N_1}, \quad z_k \in \{1, 2, \dots N_k\}, k \in \{2, 3, \dots, n\}. \notag
\end{align}

  The matrix $S_R$ contains $N_2 N_3 \cdots N_n$ blocks $B_R^{\vv{z}}$, each of which contains $N_1$ directions. Thus, $S_R$ contains $N_1N_2 \cdots N_n=N$ columns in total. Note that a block is identified using a vector $\vv{z}$ containing $n-1$ components labeled $z_2, z_3, \dots, z_n.$
Next, we provide an example in $\R^2$ to get the reader accustomed to the notation.

\begin{example}\label{ex:SR}
Select the reference point $x^0=\bbm 0&0\ebm^\top$ and the sample region $[0, 12] \times [0, 6]$.  Then $d = [12, 6]^\top = [\de_1, \de_2]^\top$. Suppose the longer side is divided into three equal parts, so $N_1=3$, and the shorter side is divided into two equal parts, so $N_2=2$.   Thus, the side lengths of one partition are $\od_1=12/3=4$ and $\od_2=6/2=3.$ The matrix $S_R$ is given by 
\begin{align*}
    S_R&=\bbm B_R^1&B_R^2\ebm,
\end{align*}
where 
\begin{align*}
    B_R^1&=\bbm \frac{12}{3}&0\\ 0&\frac{6}{2} \ebm  \bbm 1&2&3\\1&1&1 \ebm = \bbm 4&8&12\\3&3&3 \ebm, \\
    B_R^2&=\bbm \frac{12}{3}&0\\ 0&\frac{6}{2} \ebm \bbm 1&2&3\\2&2&2 \ebm =\bbm  4&8&12 \\ 6&6&6 \ebm.
\end{align*}
The sample points $x^j$ are obtained by setting  $x^j=x^0+S_Re_j$ for all $j \in \{1, 2, \dots, 6\}.$ We see that $x^1, x^2, x^3$ are associated with $B_R^1$ and $x_4, x_5, x_6$ are associated with $B_R^2.$ Figure \ref{fig:SRrectangle} illustrates the sample points built from $S_R.$

\begin{figure}[ht]
\centering
\includegraphics[width=0.5\textwidth]{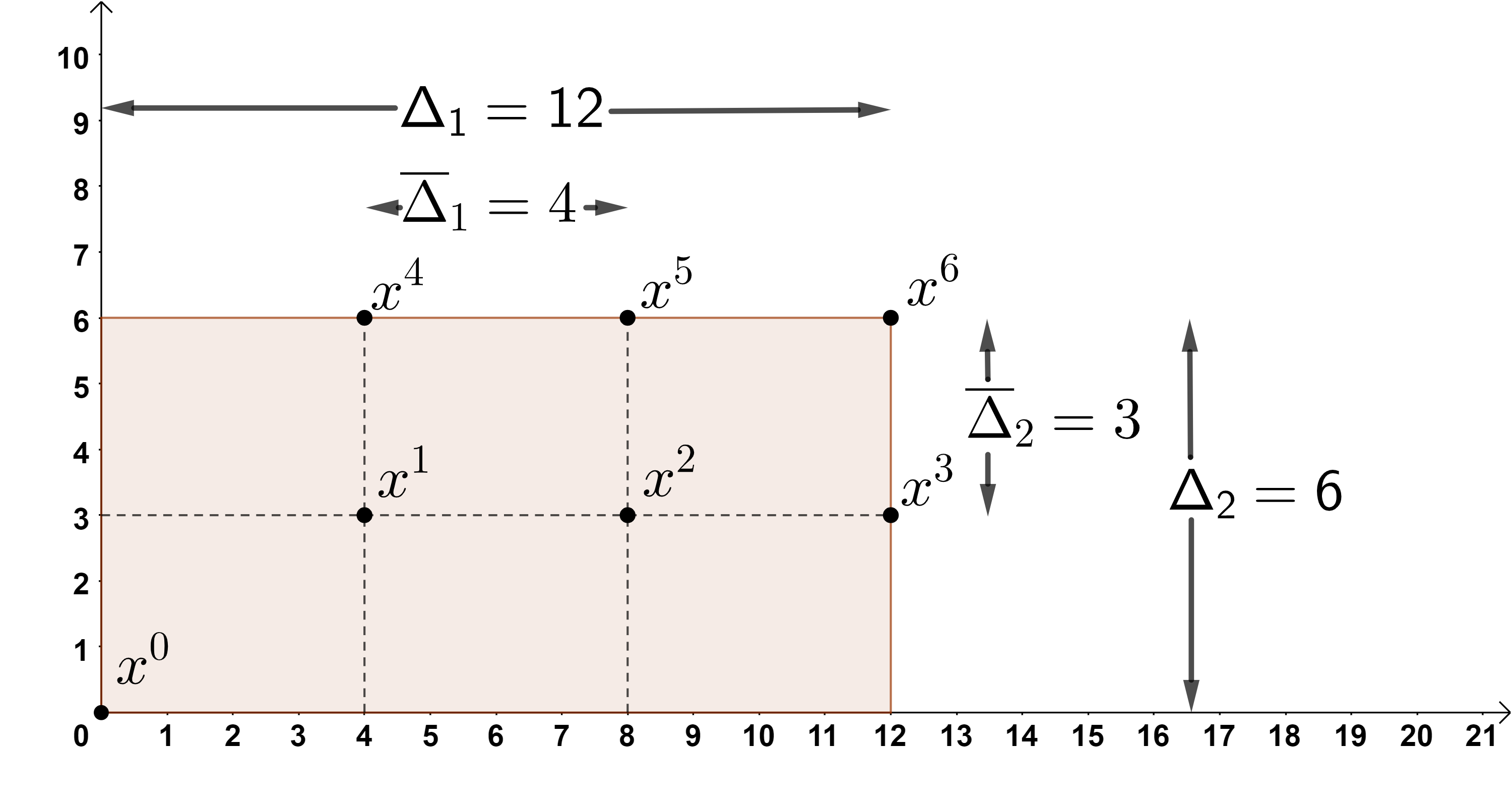} 
\caption{An example of a sample set built from a matrix $S_R$ in $\R^2$}
\label{fig:SRrectangle}
\end{figure}
\end{example}

The first goal of this paper is to find an expression for the GSG over $R(x_0;d)$ when the set of sample points forms a dense subset of $\R^n$, that is, as $N_i \rightarrow \infty$ for each $i$. As long as each $N_i$ tends to infinity (not necessarily at the same speed), we will show that an expression for the GSG over $R(x^0;d)$ can be found. Recall that the formula \eqref{eq:GSG} to compute the GSG of $f$ at $x^0$ over $S$ is 
\begin{align*}
    \nabla_s f(x^0;S)&= (S^\dagger)^\top \delta_f(x^0;S).
\end{align*}
When $S$ is full row rank, then the above formula can be written as
\begin{equation}
\nabla_s f(x^0;S) = \left (S^\top (SS^\top)^{-1} \right )^\top \delta_f(x^0;S) =(SS^\top)^{-\top}S \delta_f (x^0;S). \label{eq:Sdaggertrans}
\end{equation}
Note that $S_R \in \R^{n \times N}$ as defined in \eqref{eq:SR} is full row rank whenever $N_i\geq 2$ for all $i \in \{1, 2, \dots, n\}.$ For the remainder of this section and in Section \ref{sec:error}, we will assume $N_i\geq 2$ for all $i.$ In the following lemma, we begin our investigation of \eqref{eq:Sdaggertrans} by finding  an expression for the matrix $S_RS_R^\top.$ 

\begin{lemma}  \label{lem:sumBR}
Let $S_R \in \R^{n \times N}$ be defined as in \eqref{eq:SR}. Then
\begin{align*}
    S_R S_R^\top&=\sum_{z_2=1}^{N_2} \sum_{z_3=1}^{N_3} \cdots \sum_{z_{n}=1}^{N_n} B_R^{\vv{z}} \left (B_R^{\vv{z}} \right )^\top=U_n \in \R^{n \times n},
\end{align*}
where
\begin{align*}
[U_n]_{i,i}&=N  \frac{(N_i+1)(2N_i+1)}{6} \od_i^2, \quad i \in \{1, 2, \dots, n\}, \\
[U_n]_{i,j}&=N \frac{(N_i+1)(N_j+1)}{4} \od_i \od_j, \quad   i,j \in \{1, 2, \dots, n\}, \, i\neq j. 
\end{align*}
\end{lemma}
\begin{proof} The proof is by induction on $n$.  First, we prove the case $n=2$. We have
\begin{align*}
    \sum_{z_2=1}^{N_2} B_R^{z_2} \left ( B_R^{z_2} \right )^\top&=\bbm \od_1&0\\0&\od_2\ebm \left ( \sum_{z_2=1}^{N_2}\bbm 1&2&\cdots &N_1\\ z_2&z_2&\cdots &z_2 \ebm \bbm 1&z_2\\2&z_2\\\vdots&\vdots \\N_1&z_2 \ebm \right ) \bbm \od_1&0\\0&\od_2\ebm \\
    &=\bbm \od_1&0\\0&\od_2\ebm \left ( \sum_{z_2=1}^{N_2} \bbm \frac{N_1(N_1+1)(2N_1+1)}{6}& z_2 \frac{N_1(N_1+1)}{2}\\z_2 \frac{N_1(N_1+1)}{2}&z_2^2N_1\ebm \right ) \bbm \od_1&0\\0&\od_2\ebm\\
    &=\bbm \od_1&0\\0&\od_2\ebm \bbm \frac{N_1N_2 (N_1+1)(2N_1+1)}{6}&\frac{N_1N_2 (N_1+1)(N_2+1)}{4}\\\frac{N_1N_2 (N_1+1)(N_2+1)}{4}& \frac{N_1N_2(N_2+1)(2N_2+1)}{6} \ebm \bbm \od_1&0\\0&\od_2\ebm\\
    &=N_1N_2 \bbm \frac{(N_1+1)(2N_1+1)}{6}\od_1^2&\frac{(N_1+1)(N_2+1)}{4}\od_1 \od_2\\\frac{(N_1+1)(N_2+1)}{4}\od_1 \od_2& \frac{(N_2+1)(2N_2+1)}{6}\od_2^2 \ebm =U_2.
\end{align*}
 Now, let $n=k$ for some  $k \in\N\setminus\{1\}$. For clarity, we will write $B_R^{\vv{z_k}}$ instead of $B_R^{\vv{z}}$ to make it clear that the last component in the vector of indices $\vv{z}$ is $z_k$. Suppose  that the induction hypothesis is true for $n=k.$ We want to show that it is true for $n=k+1.$ We have
\begin{align}
    \sum_{z_2=1}^{N_2} \sum_{z_3=1}^{N_3}  \cdots  \sum_{z_{k+1}=1}^{N_{k+1}} B_R^{\vv{z_{k+1}}} \left (B_R^{\vv{z_{k+1}}} \right )^\top &=\sum_{z_2=1}^{N_2} \sum_{z_3=1}^{N_3}  \cdots  \sum_{z_{k+1}=1}^{N_{k+1}}\bbm B_R^{\vv{z_{k}}}\\ \od_{k+1} z_{k+1} \ones_{N_1}^\top \ebm \bbm \left (B_R^{\vv{z_{k}}} \right )^\top&\od_{k+1} z_{k+1} \ones_{N_1}\ebm  \notag \\
    &=\sum_{z_2=1}^{N_2} \sum_{z_3=1}^{N_3}  \cdots  \sum_{z_{k+1}=1}^{N_{k+1}} \bbm B_R^{\vv{z_{k}}} (B_R^{\vv{z_{k}}})^\top &\od_{k+1} z_{k+1} B_R^{\vv{z_{k}}} \ones_{N_1}\\ \od_{k+1} z_{k+1} \ones_{N_1}^\top \left (B_R^{\vv{z_{k}}}\right )^\top& N_1 z_{k+1}^2 \od_{k+1}^2 \ebm. \label{eq:bigsum}
\end{align}
Now, we compute the last sum of the $k$-tuple sum. We have
\begin{align}
    &\sum_{z_{k+1}=1}^{N_{k+1}} \bbm B_R^{\vv{z_{k}}} (B_R^{\vv{z_{k}}})^\top & \od_{k+1} z_{k+1} B_R^{\vv{z_{k}}} \ones_{N_1}\\ \od_{k+1} z_{k+1} \ones_{N_1}^\top \left (B_R^{\vv{z_{k}}}\right )^\top& N_1 z_{k+1}^2 \od_{k+1}^2\ebm  \notag\\
    =&\bbm N_{k+1} B_R^{\vv{z_{k}}} \left ( B_R^{\vv{z_{k}}}\right )^\top& \od_{k+1} \frac{N_{k+1} (N_{k+1}+1)}{2} B_R^{\vv{z_{k}}} \ones_{N_1} \\ \od_{k+1} \frac{N_{k+1} (N_{k+1}+1)}{2}\ones_{N_1}^\top  \left (B_R^{\vv{z_{k}}} \right )^\top& N_1 \od_{k+1}^2 \frac{N_{k+1} (N_{k+1}+1)(2N_{k+1}+1)}{6} \ebm.\label{eq:lastsum}
\end{align}
Substituting \eqref{eq:lastsum} into  \eqref{eq:bigsum}, we obtain
\begin{align} \label{eq:almostthere}
     \sum_{z_2=1}^{N_2}  \cdots \sum_{z_{k+1}=1}^{N_{k+1}} B_R^{\vv {z_{k+1}}} \left (B_R^{\vv {z_{k+1}}} \right )^\top
     &=  \sum_{z_2=1}^{N_2}  \cdots \sum_{z_{k}=1}^{N_k} \bbm N_{k+1} B_R^{\vv {z_{k}}} \left ( B_R^{\vv {z_{k}}}\right )^\top& \od_{k+1}\frac{N_{k+1}(N_{k+1}+1)}{2} B_R^{\vv {z_{k}}} \ones_{N_1} \\\od_{k+1} \frac{N_{k+1} (N_{k+1}+1)}{2} \ones_{N_1}^\top \left (B_R^{\vv {z_{k}}} \right )^\top& \od_{k+1}^2 N_1\frac{N_{k+1}(N_{k+1}+1)(2N_{k+1}+1)}{6} \ebm.
\end{align}
Let us find an expression for each of the four blocks in \eqref{eq:almostthere}. Using the induction assumption, we have 
\begin{align*}
    \sum_{z_2=1}^{N_2} \sum_{z_3=1}^{N_3} \cdots \sum_{z_k}^{N_k} N_{k+1}B_R^{\vv{z_k}}\left ( B_R^{\vv{z_k}}\right )=N_{k+1}U_k.
\end{align*}
The off-diagonal blocks are transposes of each other, with 
\begin{align}
    \sum_{z_2=1}^{N_2} \sum_{z_3=1}^{N_3} \cdots \sum_{z_k=1}^{N_k} \od_{k+1}\frac{N_{k+1}(N_{k+1}+1)}{2} B_R^{\vv {z_{k}}} \ones_{N_1}&=\od_{k+1}\frac{N_{k+1}(N_{k+1}+1)}{2}\bbm \od_1 N_1 \frac{(N_1+1)}{2} N_2 N_3 \cdots N_k\\ \od_2 N_1 N_2 \frac{(N_2+1)}{2} N_3 N_4 \cdots N_k \\ \vdots \\ \od_k N_1 N_2 \cdots N_{k-1} N_k \frac{(N_k+1)}{2}\ebm \notag \\
    &=\bbm N \od_1 \od_{k+1} \frac{(N_1+1)(N_{k+1}+1)}{4}\\  N \od_2 \od_{k+1} \frac{(N_2+1)(N_{k+1}+1)}{4} \\\vdots \\  N \od_k \od_{k+1} \frac{(N_k+1)(N_{k+1}+1)}{4}\ebm.\notag 
\end{align}
The last block is
\begin{align*}
    \sum_{z_2=1}^{N_2}\sum_{z_3=1}^{N_3} \cdots \sum_{z_k=1}^{N_k} \od_{k+1}^2 N_1\frac{N_{k+1}(N_{k+1}+1)(2N_{k+1}+1)}{6}&=N \od_{k+1}^2\frac{(N_{k+1}+1)(2N_{k+1}+1)}{6}.
\end{align*}
Hence, \eqref{eq:almostthere} is equal to $U_{k+1}.$
By the principle of induction, the claim is true for all $n \in\N\setminus\{1\}$.
\end{proof}

The GSG uses the inverse of $(S_RS_R^\top)^\top$.  The following lemma finds an expression for this inverse. 

\begin{lemma} \label{lem:invSRSRt} 

Let $S_R \in \R^{n \times N}$ be defined as in \eqref{eq:SR}. Then 
\begin{align*}
   (S_RS_R^\top)^{-\top}&= \frac{12}{N} \left ( E -\frac{3}{1+3s} yy^\top \right ),
\end{align*}
where $$E= \diag \bbm \frac{1}{(N_1^2-1) \od_1^2}& \cdots&  \frac{1}{(N_n^2-1)\od_n^2} \ebm \in \R^{n \times n}, ~s=\sum_{i=1}^n \frac{N_i+1}{N_i-1}, ~\text{and} \quad y=\bbm \frac{1}{(N_1-1)\od_1} &\cdots& \frac{1}{(N_n-1)\od_n}\ebm^\top \in \R^n.$$
\end{lemma}
\begin{proof}
Note that for all $n \in\N \setminus\{1\}$ the symmetric matrix $U_n$  can be written as
\begin{align*}
    U_n&=N \left ( \diag \bbm \frac{(N_1-1)(N_1+1)}{12} \od_1^2& \cdots&\frac{(N_n-1)(N_n+1)}{12} \od_n^2\ebm    + \bbm \frac{(N_1+1)}{2}\od_1\\ \vdots \\ \frac{(N_n+1)}{2}\od_n \ebm \bbm \frac{(N_1+1)}{2}\od_1& \cdots& \frac{(N_n+1)}{2}\od_n \ebm \right ).
\end{align*}
Let $\dot{d}=\bbm \frac{(N_1+1)}{2}\od_1& \cdots& \frac{(N_n+1)}{2}\od_n \ebm^\top \in \R^n$ and $\widetilde{D}=\diag \bbm \frac{N_1^2-1}{12} \od_1^2& \cdots&\frac{N_n^2-1}{12} \od_n^2\ebm \in \R^{n \times n}.$ Then 
\begin{align*}
    U_n=N(\widetilde{D}+\dot{d} \, \dot{d}^\top).
\end{align*}
Therefore, using \eqref{eq:sherman}, we obtain
\begin{align}
    (U_n)^{-\top}&=(U_n^\top )^{-1}=U_n^{-1} \notag \\
    &=\frac{1}{N}\left (\widetilde{D}^{-1}- \frac{\widetilde{D}^{-1}\dot{d} \, \dot{d}^\top \widetilde{D}^{-1}}{1+\dot{d}^\top \widetilde{D}^{-1} \dot{d}} \right ). \label{eq:uninv}
\end{align}
 The denominator of the second term in \eqref{eq:uninv} is
 \begin{align*}
     1+\dot{d}^\top \widetilde{D}^{-1} \dot{d}&= 1+\bbm \frac{(N_1+1)}{2}\od_1&\cdots&\frac{(N_n+1)}{2}\od_n \ebm \diag \bbm \frac{12}{(N_1^2-1)\od_1^2}&\cdots&\frac{12}{(N_n^2-1)\od_n^2} \ebm \bbm \frac{(N_1+1)}{2}\od_1\\ \vdots \\ \frac{(N_n+1)}{2} \od_n \ebm \\
     &=1+ \bbm \frac{6}{(N_1-1)\od_1}&\cdots&\frac{6}{(N_n-1)\od_n}\ebm \bbm \frac{(N_1+1)}{2}\od_1\\ \vdots \\ \frac{(N_n+1)}{2} \od_n \ebm \\
     &=1+3 \sum_{i=1}^n\frac{N_i+1}{N_i-1}=1+3s.
 \end{align*}
The numerator of the second term in \eqref{eq:uninv} is
 \begin{align*}
     \widetilde{D}^{-1}\dot{d} \, \dot{d}^\top \widetilde{D}^{-1}&=\bbm \frac{6}{(N_1-1)\od_1}\\ \vdots \\ \frac{6}{(N_n-1)\od_n}\ebm  \bbm \frac{6}{(N_1-1)\od_1}& \cdots & \frac{6}{(N_n-1)\od_n}\ebm= 36yy^\top.
 \end{align*}
 Since $\widetilde{D}^{-1}=12E,$ we get
 \begin{align*}
     (U_n)^\top&= \frac{1}{N}\left ( \widetilde{D}^{-1} - \frac{36}{1+3s}yy^\top\right)=\frac{12}{N} \left (E-\frac{3}{1+3s}yy^\top \right ).\qedhere
 \end{align*}
\end{proof}

Using the previous lemma, we can now provide an expression for the transpose of the Moore--Penrose inverse $S_R^\dagger$.

\begin{corollary}[The matrix $(S_R^\dagger)^\top$] \label{prop:SRdt}
Let $S_R \in \R^{n \times N}$ be defined as in \eqref{eq:SR}. Then 
\begin{align*}
    (S_R^\dagger)^\top&=\frac{12}{N} \left (E-\frac{3}{1+3s}yy^\top \right ) S_R,
\end{align*}
where $E=\diag \bbm \frac{1}{(N_1^2-1)\od_1^2}&\cdots&\frac{1}{(N_n^2-1)\od_n^2}\ebm \in \R^{n \times n}$, $y=\bbm \frac{1}{(N_1-1)\od_1}&\cdots& \frac{1}{(N_n-1)\od_n}\ebm^\top \in \R^n,$ and the scalar $s=\sum_{i=1}^n \frac{N_i+1}{N_i-1}.$
\end{corollary}
\begin{proof}
Since $S_R$ has full row rank, we have
\begin{align*}
    (S_R^\dagger)^\top&=\left ( S_R^\top (S_R S_R^\top)^{-1}  \right )^\top=(S_R S_R)^{-\top}S_R=\frac{12}{N} \left (E-\frac{3}{1+3s}yy^\top \right ) S_R,
\end{align*}
by Lemma \ref{lem:invSRSRt}.
\end{proof}

In the following proposition, we investigate the limit of $N (S_R S_R^\top)^{-\top}$ when all $N_i$ go to infinity. Applying these limits, the sample points contained in $R(x_0;d)$ form a dense grid of $\R^n.$ To make the notation compact, we write $\lim_{\vv{N} \to \vv{\infty}}$ to represent the limit as $N_1 \to \infty, N_2 \to \infty, \dots, N_n \to \infty.$ The reason we are interested in the limit of $N (S_RS_R^\top)^{-\top}$ is that, assuming the limits exist, we  may write
\begin{align}
    \limnvec \nabla_s f(x^0;S_R)&=\limnvec (S_R^\dagger)^\top \delta_f (x^0;S_R) \notag\\
    &=\limnvec \frac{N}{\de}(S_RS_R^\top)^{-\top} \frac{\de}{N} S_R \delta_f (x^0;S_R) \notag\\
    &=\frac{1}{\de}\left [ \limnvec  N (S_RS_R^\top)^{-\top} \right ] \left [\limnvec \frac{\de}{N} S_R \delta_f (x^0;S_R) \right ], \label{eq:twolimits}
\end{align}
where $\de=\de_1\de_2\cdots\de_n.$ So if we can show that the two limits in \eqref{eq:twolimits} exist, then we have found the limit of the GSG over a dense hyperrectangle. Note that the term $\Delta$ remains inside the second limit, as we will show that the expression in the second limit is a $n$-tuple Riemann sum. 
\begin{proposition} \label{prop:limuninv}
Let $S_R$ be defined as in \eqref{eq:SR}. Then 
\begin{align*}
    \lim_{\vv{N} \to \vv{\infty}} N (S_RS_R^\top)^{-\top}&=L_n\in \R^{n \times n} 
\end{align*}
where  
\begin{align*}
    [L_n]_{i,i}&=\frac{12(3n-2)}{\de_i^{2}(3n+1)}, \quad i \in \{1, 2, \dots, n\}, \\
    [L_n]_{i,j}&= \frac{-36}{\de_i \de_j (3n+1)}, \quad i,j \in \{1, 2, \dots, n\}, i \neq j.
\end{align*}
\end{proposition}
\begin{proof}
By Corollary \ref{prop:SRdt} we have
\begin{align}
    \lim_{\vv{N} \to \vv{\infty}} N (S_RS_R^\top)^{-\top}&=\lim_{\vv{N} \to \vv{\infty}} N \frac{12}{N} \left ( E -\frac{3}{1+3s}y y^\top \right ) =\lim_{\vv{N} \to \vv{\infty}} 12 \left ( E -\frac{3}{1+3s}y y^\top \right ) \label{eq:limsrsr},
\end{align}
where $E=\diag \bbm \frac{1}{(N_1^2-1)\od_1^2}&\cdots& \frac{1}{(N_n^2-1)\od_n^2}\ebm \in \R^{n \times n},$ $y=\bbm \frac{1}{(N_1-1)\od_1}&\cdots&\frac{1}{(N_n-1)\od_n}\ebm^\top \in \R^n$ and the scalar  $s=\sum_{i=1}^{n} \frac{N_i+1}{N_i-1}.$ We show that this converges componentwise to $L$. 

We begin with the diagonal entries. 
Applying $\od_i^2 =N_i^2/\de_i^2$, note that
\begin{align*}
    [yy^\top]_{i,i}&=\frac{1}{(N_i-1)^2\od_i^2}=\frac{N_i^2}{(N_i-1)^2\de_i^2}.
\end{align*}
Substituting this and $\od_i^2=N_i^2/\de_i^2$ into the definition of $E$ yields
\begin{align*}
     \lim_{\vv{N}\to \vv{\infty}}\left[12\left (E-\frac{3}{1+3s}y y^\top \right ) \right ]_{i,i}&=12 \limnvec \left( \frac{N_i^2}{(N_i-1)(N_i+1)\de_i^2}-\frac{3}{1+3\sum_{i=1}^n \frac{N_i+1}{N_i-1}}\frac{N_i^2}{(N_i-1)^2\de_i^2}\right) \\
     &=12 \left ( \frac{1}{\de_i^2}-\frac{3}{(1+3n)\de_i^2}\right )\\
     &=\frac{12(3n-2)}{\de_i^2(3n+1)}, \quad i \in \{1, 2, \dots, n\}.
\end{align*}
Similarly, the off-diagonal entries of the matrix in \eqref{eq:limsrsr} are given by
\begin{align*}
    \limnvec \left [\frac{-36}{1+3\sum_{i=1}^n \frac{N_i+1}{N_i-1}} \frac{N_i N_j}{(N_i-1)(N_j-1)\de_i \de_j}\right ]_{i,j}&=\frac{-36}{(1+3n)\de_i \de_j}, \quad i,j \in \{1, 2, \dots, n\}, i\neq j.\qedhere
\end{align*}
\end{proof}

\subsection{Generalization: using an arbitrary point in each partition}

We now generalize $S_R$ to a matrix that allows choosing an arbitrary point in each partition of the sample region, not necessarily the right endpoint of each partition.  The matrix containing  all directions used to obtain an arbitrary sample point in each partition will be denoted by $S.$ Let $N=N_1N_2 \cdots N_n$ and $\overline{D}=\diag \bbm \od_1&\od_2&\cdots&\od_n\ebm \in \R^{n \times n}$ where $\od_i=\de_i/N_i$ for all $i$. The matrix $S \in \R^{n \times N}$ can be written as a block matrix in the following way:
\begin{align}\label{eq:SA}
    S&=\bbm B^{1, 1,\dots,1,1}&B^{1, 1, \dots,1,2}&\cdots&B^{N_2,N_3, \dots, N_{n-1}, N_n} \ebm
\end{align}
where
\begin{align*}
    B^{\vv{z}}&=B_R^{\vv{z}}-B_M^{\vv{z}} \in \R^{n \times N_1}.
\end{align*}
The block $B_R^{\vv{z}} \in \R^{n \times N_1}$ is defined in \eqref{eq:BRabm} and the block $B_M^{\vv{z}} \in \R^{n \times N_1}$ is
\begin{align*}
    B_M^{\vv{z}}&=\overline{D} \ddot{B}_M^{\vv{z}}
\end{align*}
where  all entries of $\ddot{B}_M^{\vv{z}}$ are in $[0, 1].$ 
Let $S_M \in \R^{n \times N}$ be defined as 
\begin{align*}
    S_M&= \bbm B_M^{1, 1, \dots, 1, 1} & B_M^{1, 1, \dots, 1, 2}&\cdots& B_M^{N_2, N_3, \dots, N_{n-1}, N_n}\ebm.
\end{align*}
Then $S$ can be written as $S=S_R-S_M.$ 
Let us provide an example of $S$ in $\R^2$. 

\begin{example}
Consider the same sample region as Example \ref{ex:SR}.   Select the `arbitrary' points $[2, 2]^\top,$ $ [5, 1]^\top,$ $ [8, 3]^\top,$ $ [1, 1]^\top,$ $ [6, 4.5]^\top,$ and $[12, 6]^\top$.  These are visualized in Figure \ref{fig:Srandom}.

\begin{figure}[ht]
\centering
\includegraphics[width=0.6\textwidth]{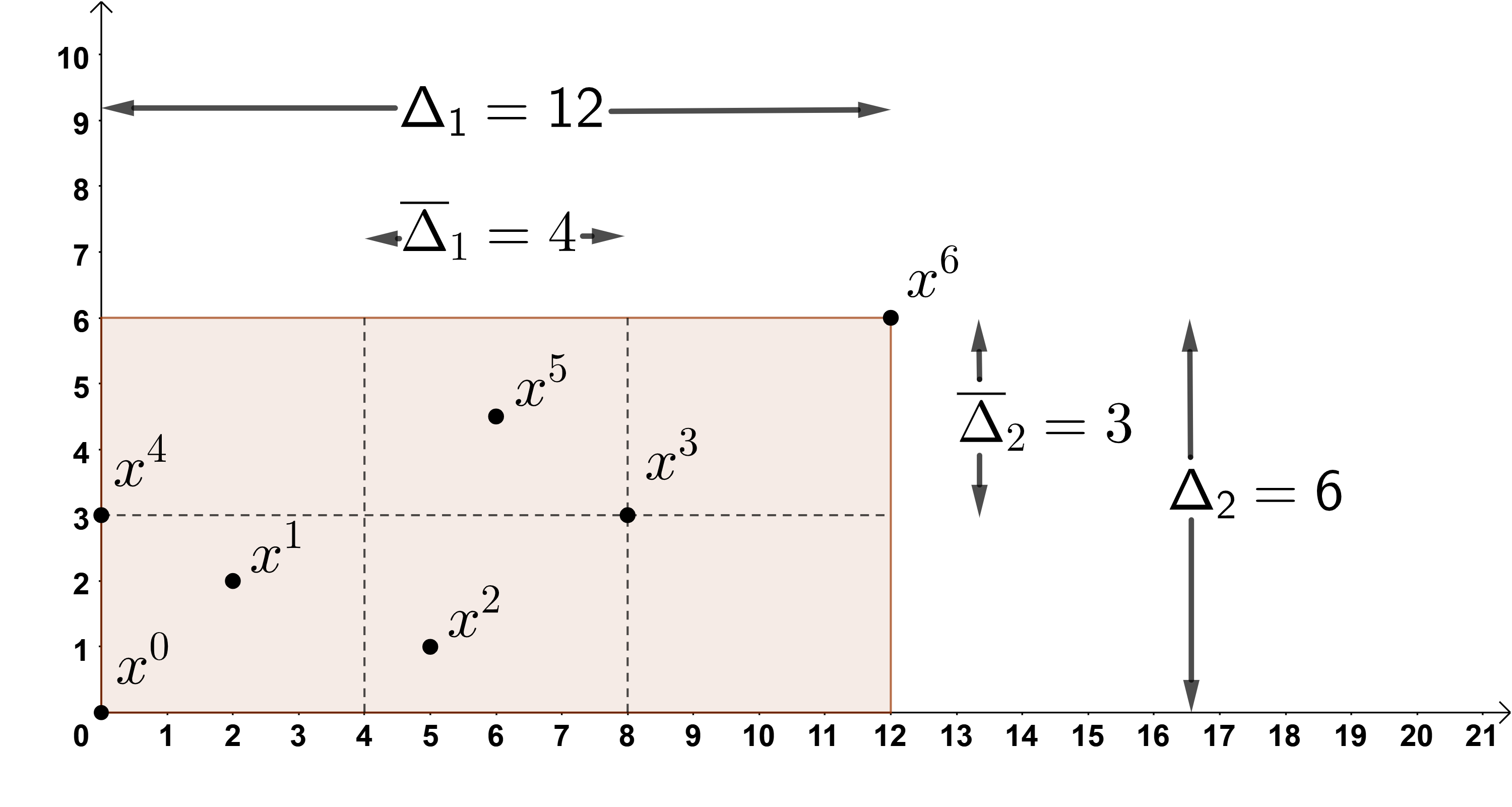} 
\caption{An example of sample set built from a matrix $S$ in $\R^2$}
\label{fig:Srandom}
\end{figure}

The matrix $S \in \R^{2 \times 6}$ is  given by
\begin{align*}
    S&=\bbm B^1&B^2\ebm,
\end{align*}
where
\begin{align*}
    B^1&=B_R^1-B_M^1\\
    &=\bbm \frac{12}{3}&0\\0&\frac{6}{2}\ebm \bbm 1&2&3\\1&1&1 \ebm-\bbm \frac{12}{3}&0\\0&\frac{6}{2}\ebm \bbm \frac{1}{2}&\frac{3}{4}&1\\\frac{1}{3}&\frac{2}{3}&0 \ebm \\
    &=\bbm4&8&12\\3&3&3 \ebm-\bbm2&3&4\\1&2&0 \ebm =\bbm 2&5&8\\2&1&3\ebm,
\end{align*}
and
\begin{align*}
     B^2&=B_R^2-B_M^2\\
    &=\bbm \frac{12}{3}&0\\0&\frac{6}{2}\ebm \bbm 1&2&3\\2&2&2 \ebm-\bbm \frac{12}{3}&0\\0&\frac{6}{2}\ebm \bbm 1&\frac{1}{2}&0\\1&\frac{1}{2}&0 \ebm \\
    &=\bbm4&8&12\\6&6&6 \ebm-\bbm4&2&0\\3&\frac{3}{2}&0 \ebm =\bbm 0&6&12\\3&\frac{9}{2}&6\ebm.
\end{align*}
The sample points $x^j$ are built by setting $x^j=x^0+Se_j$ for all $j \in \{1, 2, \dots, 6\}.$ We see that $x^1, x^2, x^3$ are associated with $B^1$ and $x^4, x^5, x^6$ are associated with $B^2.$ 
\end{example}

The following proposition generalizes Proposition \ref{prop:limuninv} by considering $S$ in place of $S_R.$
\begin{proposition}  \label{prop:limSASAti}
Let  $S \in \R^{n \times N}$ be defined as in  \eqref{eq:SA}. Then 
\begin{align*}
    \limnvec N (S S^\top)^{-\top}=L_n \in \R^{n \times n},
\end{align*}
where $L_n$ is defined as in Proposition \ref{prop:limuninv}.
\end{proposition}
\begin{proof}
We have
\begin{align*}
    \limnvec N(S S^\top)^{-\top}&=  \limnvec \left (\frac{1}{N} S S^\top \right )^{-\top}.
\end{align*}
Note that the inverse of $S S^\top$ is well-defined, since $S$ is full row rank  whenever $N_i\in\N\setminus\{1\}$ for all $i.$ It follows that $SS^\top$ is full rank, so the inverse of $(SS^\top)^\top$ exists.  Since the inverse of  $(S S^\top)^\top$ is a continuous function  with respect to $\vv{N},$ we may take the limit inside the inverse. We obtain 
\begin{align}
      \limnvec \left (\frac{1}{N} S S^\top \right )^{-\top}&= \left ( \limnvec \frac{1}{N} S S^\top\right )^{-\top} \notag \\
      &= \left ( \limnvec \frac{1}{N}S_R S_R^\top -\limnvec \frac{1}{N} S_M S_R^\top - \limnvec \frac{1}{N} S_RS_M^\top+ \limnvec \frac{1}{N} S_MS_M^\top \right )^{-\top}. \label{eq:limexpanded}
\end{align}
Now, we show that $\limnvec \frac{1}{N} S_M S_R^\top, \limnvec \frac{1}{N} S_RS_M^\top,$ and $\limnvec \frac{1}{N} S_MS_M^\top$ are equal to the $n \times n$ zero matrix.

We begin with showing $\limnvec \frac{1}{N} S_MS_M^\top=0_{n \times n}$.  Let $\overline{D} \in \R^{n\times n}$ be the diagonal matrix with entries $\de_1/N_1, \dots, \de_n/N_n.$  We have 
\begin{align}
    \frac{1}{N}S_MS_M^\top&=\frac{1}{N}\sum_{z_2=1}^{N_2} \cdots \sum_{z_n=1}^{N_n} B_M^{\vv{z}}B_M^{\vv{z}} \notag\\
    &=\frac{1}{N} \overline{D} \left (\sum_{z_2=1}^{N_2}\cdots \sum_{z_n=1}^{N_n} \ddot{B}_M^{\vv{z}} (\ddot{B}_M^{\vv{z}})^\top \right )\overline{D}.\label{eq:SMSM}
\end{align}
Since all entries in the matrix $\ddot{B}_M^{\vv{z}}$ are  contained in $[0,1],$ the ($n-1$)-tuple sum in \eqref{eq:SMSM} is bounded componentwise below by the matrix $0_{n \times n}$ and above by
\begin{align*}
    \sum_{z_2=1}^{N_2}\cdots \sum_{z_n=1}^{N_n} \ddot{B}_M^{\vv{z}} (\ddot{B}_M^{\vv{z}})^\top&\leq  \sum_{z_2=1}^{N_2}\cdots \sum_{z_n=1}^{N_n}  \ones_n \ones_{N_1}^\top \ones_{N_1} \ones_n^\top \\
    &=N_1 \sum_{z_2=1}^{N_2}\cdots \sum_{z_n=1}^{N_n} \ones_n \ones_n^\top = N \ones_n \ones_n^\top.
\end{align*}
It follows that componentwise
\begin{align*}
    \limnvec \frac{1}{N} \overline{D} \left (\sum_{z_2=1}^{N_2}\cdots \sum_{z_n=1}^{N_n} \ddot{B}_M^{\vv{z}} (\ddot{B}_M^{\vv{z}})^\top \right )\overline{D} &\leq \limnvec \frac{1}{N} \overline{D} N \ones_n \ones_n^\top \overline{D}\\
    &=\limnvec \bbm \frac{\de_1}{N_1}\\\frac{\de_2}{N_2}\\\vdots\\\frac{\de_n}{N_n}\ebm  \bbm \frac{\de_1}{N_1}&\frac{\de_2}{N_2}&\cdots&\frac{\de_n}{N_n}\ebm=0_{n \times n}.
\end{align*}
By the Squeeze Theorem, $\limnvec \frac{1}{N}S_MS_M^\top=0_{n \times n}.$ 

Now, we show that $\limnvec \frac{1}{N} S_RS_M^\top=0_{n \times n}.$ We have
\begin{align}\label{eq:BRBM}
    \frac{1}{N}S_RS_M^\top&=\frac{1}{N} \overline{D} \left (\sum_{z_2=1}^{N_2} \cdots \sum_{z_n=1}^{N_n} \ddot{B}_R^{\vv{z}}(\ddot{B}_M^{\vv{z}})^\top\right )\overline{D}.
\end{align}
The ($n-1$)-tuple sum in \eqref{eq:BRBM} is bounded componentwise below by $0_{n \times n}$.  A componentwise upper bound for \eqref{eq:BRBM} is 
\begin{align*}
    \frac{1}{N} \overline{D} \left (\sum_{z_2=1}^{N_2} \cdots \sum_{z_n=1}^{N_n} \ddot{B}_R^{\vv{z}}(\ddot{B}_M^{\vv{z}})^\top\right )\overline{D} &\leq \frac{1}{N} \overline{D} \left ( \sum_{z_2=1}^{N_2} \cdots \sum_{z_n=1}^{N_n} \ddot{B}_R^{\vv{z}} \ones_n \ones_{N-1}^\top \right ) \overline{D} \\
    &=\frac{1}{N}  \overline{D} \frac{N}{2} \bbm N_1+1&N_1+1&\cdots& N_1+1\\ N_2+1&N_2+1&\cdots&N_2+1\\ \vdots&\vdots&\ddots&\vdots\\N_n+1&N_n+1&\cdots&N_n+1\ebm \overline{D} \\
    &=\frac{1}{2} \bbm \frac{\de^2_1(N_1+1)}{N_1^2}&\frac{\de_1 \de_2 (N_1+1)}{N_1 N_2}&\cdots& \frac{\de_1 \de_n (N_1+1)}{N_1 N_n}\\ \frac{\de_1 \de_2(N_2+1)}{N_1 N_2}&\frac{\de^2_2 (N_2+1)}{N_2^2}&\cdots&\frac{\de_2 \de_n(N_2+1)}{N_2 N_n}\\ \vdots&\vdots&\ddots&\vdots\\\frac{\de_1 \de_n(N_n+1)}{N_1N_n}&\frac{\de_2 \de_n (N_n+1)}{N_2 N_n}&\cdots&\frac{\de_n^2(N_n+1)}{N_n^2}\ebm.
\end{align*}
It follows that $\limnvec S_RS_M^\top \leq 0_{n \times n},$ and by the Squeeze Theorem, $\limnvec S_RS_M^\top = 0_{n \times n}.$ As $S_MS_R^\top = (S_RS_M^\top)^\top$, we also have $\limnvec S_MS_R^\top=0_{n \times n}.$ Thus, \eqref{eq:limexpanded} reduces to
\begin{align*}
    \left ( \limnvec \frac{1}{N}S_R S_R^\top -\limnvec \frac{1}{N} S_M S_R^\top - \limnvec \frac{1}{N} S_RS_M^\top+ \limnvec \frac{1}{N} S_MS_M^\top \right )^{-\top}&=\left ( \limnvec \frac{1}{N} S_RS_R^\top \right )^{-\top}\\
    &=\limnvec N \left ( S_RS_R^\top\right )^{-\top}=L_n
\end{align*}
by Proposition  \ref{prop:limuninv}.
\end{proof}


The previous proposition gives an expression for the first limit in \eqref{eq:twolimits}. The following theorem gives the limit of the product $\frac{\de}{N} S\delta_f(x^0;S)$, the second limit in \eqref{eq:twolimits}, as a multiple integral over $R(0;d)$.

\begin{proposition} \label{prop:tn}
Let $f \in \mathcal{C}^0$ on an open domain containing  $R(x^0; d) \subseteq \dom f.$ Let $d=\bbm \de_1&\de_2&\cdots&\de_n \ebm^\top>0,$ $x=\bbm x_1&\cdots&x_n \ebm^\top$ and    let $S \in \R^{n \times N}$ be defined as in \eqref{eq:SA}. Then 
\begin{align*}
    \limnvec \frac{\de}{N} S \, \ds (x^0;S)&= T_n \in \R^n,
\end{align*}
where
\begin{align*}
    [T_n]_i&= \nintrecc x_i \left ( f(x^0+x)-f(x^0)\right )\, dx, \quad i \in \{1, 2, \dots, n\}.
\end{align*}
\end{proposition}
\begin{proof}
We have 
\begin{align}
    \limnvec \frac{\de}{N} S \, \ds (x^0;S)&= \limnvec \frac{\de}{N} \bbm B^{1, 1, \dots, 1}&\cdots B^{N_2, N_3, \dots, N_n} \ebm  \bbm \ds (x^0;B^{1, 1, \dots, 1}) \\ \vdots \\ \ds (x^0;B^{N_2, N_3, \dots, N_n}) \ebm \notag \\
    &= \limnvec   \sum_{z_2=1}^{N_2} \sum_{z_3=1}^{N_3} \cdots \sum_{z_{n}=1}^{N_n} \frac{\de}{N} B^{\vv{z}} \ds (x^0;B^{\vv{z}}). \label{eq:intform}
\end{align}
Note that $\frac{\de}{N}$ is the volume of one partition of $R(x^0;d).$ Recall that
\begin{align*}
    \delta_f(x^0;B^{\vv{z}})&=\bbm f(x^0+B^{\vv{z}}e_1)-f(x^0)&\cdots& f(x^0+B^{\vv{z}}e_{N_1})-f(x^0)\ebm^\top \in \R^{N_1}.
\end{align*}
The matrix $B^{\vv{z}}$ has dimension $n \times N_1$ so \eqref{eq:intform} is the  limit of a vector in $\R^n.$ Since $f\in \mathcal{C}^0$ on an open domain containing $R(x^0;d),$ \eqref{eq:intform} is a vector of $n$ definite integrals:
\begin{align*}
\limnvec   \sum_{z_2=1}^{N_2} \sum_{z_3=1}^{N_3} \cdots \sum_{z_{n}=1}^{N_n} \frac{\de}{N} B^{\vv{z}} \ds (x^0;B^{\vv{z}})&=\bbm \nintrecc x_1 \left ( f(x^0+x)-f(x^0) \right ) \, dx \\\nintrecc x_2 \left ( f(x^0+x)-f(x^0) \right ) \, dx  \\ \vdots \\ \nintrecc x_n \left ( f(x^0+x)-f(x^0) \right )\, dx \ebm=T_n.\qedhere
\end{align*}
\end{proof}


Now we are ready for our first main result: the limiting behaviour of the GSG at $x^0$ over the dense grid $R(x^0;d)$. The result of Theorem \ref{thm:limsimplexgrad} below is expressed as a multiple definite integral.

\begin{theorem}[Limiting behaviour of the GSG of $f$ at $x^0$ over $S$] \label{thm:limsimplexgrad}
 Let $f \in \mathcal{C}^0$ on an open domain containing $R(x^0;d) \subseteq \dom f,$ $d=\bbm \de_1&\de_2&\cdots&\de_n\ebm^\top >0.$  Let $S \in \R^{n \times N}$ be defined as in \eqref{eq:SA}. Then 
\begin{align}
    \limnvec \nabla_s f(x^0;S)&=\de^{-1} L_n T_n, \label{eq:limgsg}
\end{align}
where the entries of $L_n \in \R^{n \times n}$ are given by
\begin{align*}
    [L_n]_{i,i}&=\frac{12(3n-2)}{\de_i^2(3n+1)}, \quad i \in \{1, 2, \dots, n\}, \\
    [L_n]_{i,j}&= \frac{-36}{\de_i \de_j(3n+1)}, \quad i,j \in \{1, 2, \dots, n\}, i \neq j,
\end{align*}
and the entries of $T_n \in \R^n$ are given by 
\begin{align}
    [T_n]_i&= \nintrecc x_i \left ( f(x^0+x)-f(x^0)\right ) \, dx, \quad i \in \{1, 2, \dots, n\}. \label{eq:Tn}
\end{align}
\end{theorem}
\begin{proof}
We have
\begin{align*}
    \limnvec \nabla_s f(x^0;S)&= \limnvec (S^\dagger)^\top \ds (x^0;S) \\
    &=\limnvec \left ( S S^\top \right )^{-\top} S \ds (x^0;S) \\
    &=\limnvec \frac{N}{\de} \left ( S S^\top\right )^{-\top} \frac{\de}{N}  S \ds (x^0;S).
\end{align*}
By Proposition \ref{prop:limSASAti},
\begin{align*}
     \limnvec N(S S^\top)^{-\top}=L_n \in \R^{n \times n}.
\end{align*}
By Proposition  \ref{prop:tn},
\begin{align*}
    \limnvec \frac{\de}{N} S \ds (x^0;S)&= T_n \in \R^n.
\end{align*}
Therefore, 
\begin{align*}
    \limnvec \nabla_s f(x^0;S)&= \frac{1}{\de}\limnvec N(S S^\top)^{-\top} \limnvec \frac{\de}{N} S \ds (x^0;S) \\
    &= \de^{-1} L_n T_n.\qedhere
\end{align*}
\end{proof}

\begin{remark}
The previous result also agrees with what was found in \cite{MR4163088} for $n=1$. Indeed, in $\R$ \cite[Theorem 4.1]{MR4163088} found that 
\begin{align*}
    \limn \nabla_s f(x^0;S)=\frac{3}{\de^3} \int_{0}^{\de} x (f(x^0+x)-f(x^0)) dx,
\end{align*}
which is what \eqref{eq:limgsg} becomes when $n=1$.
\end{remark}

We end this section with an example of $\limnvec \nabla_s f(x^0;S)$.

\begin{example}
Let $f:\R^2 \to \R: x \mapsto x_1^2+x_2^2$ and the reference point $x^0=\bbm 3&1 \ebm^\top$.
Consider the sample region $R(x^0;\bbm 1&1\ebm^\top ),$ a square of side length 1.  By Theorem \ref{thm:limsimplexgrad}, we know that
\begin{align*}
    \limnvec \nabla_S f(x^0;S)&=\frac{1}{\Delta}L_2 T_2 \\
    &=\frac{1}{(1)(1)} \bbm \frac{48}{7}&\frac{-36}{7}\\\frac{-36}{7}&\frac{48}{7} \ebm \bbm \int_{0}^{1} \int_{0}^1 x_1 \left ( (3+x_1)^2+(1+x_2)^2-10 \right ) dx_1 dx_2 \\ \int_{0}^{1} \int_{0}^1 x_2 \left ( (3+x_1)^2+(1+x_2)^2-10 \right ) dx_1 dx_2 \ebm \\
    &= \bbm \frac{48}{7}&\frac{-36}{7}\\\frac{-36}{7}&\frac{48}{7} \ebm \bbm \frac{35}{12} \\ \frac{31}{12} \ebm \approx \bbm 6.71\\2.71\ebm .
\end{align*}
The absolute error is approximately $ \left \Vert \bbm 6.71&2.71\ebm^\top -\bbm6&2\ebm^\top \right  \Vert \approx 1.01.$
\end{example}

Now that we have an expression for the GSG as the number of points tends to infinity,  the next step is to define an error bound \emph{ad infinitum}. This is the focus of the next section.

\section{Error bound ad infinitum of GSG over a hyperrectangle}\label{sec:error}
In this section, we use the results obtained thus far to formulate an error bound ad infinitum that does not depend on the number of points used in the hyperrectangle $R(x^0;d)$.

To obtain an error bound ad infinitum for the GSG at $x^0$ over $R(x^0;d)$,   we require  the function  $f$ to be $\mathcal{C}^2$ on an open domain containing $R(x^0;d).$ We will write $f(x^0+x)$ as the first-order Taylor expansion. That is 
 \begin{align}
     f(x^0+x)&=f(x^0)+\nabla f(x^0)^\top x +R_1(x^0+x),~\mbox{where}~
 R_1(x^0+x)= \frac{1}{2}x^\top \nabla^2 f(\xi) x, \label{eq:firstordertaylor}
 \end{align}
 for $\xi \in \R^n$ between $x^0$ and $x^0+x$. 
 
By rewriting $f(x^0+x)$  in the form of \eqref{eq:firstordertaylor},  the components of the vector $T_n$ defined in \eqref{eq:Tn} can be written as
\begin{align*}
    [T_n]_i&=\nintrecc x_i \left (f(x^0+x)-f(x^0)\right) \dexxe \\
    &=\nintrecc x_i \left ( \nabla f(x^0)^\top x+R_1(x^0+x) \right ) \dexxe \\
    &=\nintrecc x_i \nabla f(x^0)^\top x \dexxe + \nintrecc x_i R_1(x^0+x) \dexxe.
\end{align*}
for $i \in \{1, 2, \dots, n\}.$
Let $v \in \R^n$ be the vector defined by
\begin{align}
    v_i&=\nintrecc x_i \nabla f(x^0)^\top x\dexxe, \quad i \in \{1, 2, \dots, n\}, \label{eq:vectorv}
\end{align}
and $w \in \R^n$ be the vector defined by
\begin{align*}
    w_i&=\nintrecc x_i R_1(x^0+x) \dexxe, \quad i \in \{1, 2, \dots, n\}. 
\end{align*}
Then the expression for $\displaystyle\limnvec \nabla_s f(x^0;S)$ given in \eqref{eq:limgsg} can be expressed as
\begin{equation}
    \limnvec \nabla_s f(x^0;S)
    = \frac{1}{\de} L_n v+ \frac{1}{\de} L_nw. \label{eq:limgsg2terms}
\end{equation}
We begin by showing that the first term in \eqref{eq:limgsg2terms} is equal to $\nabla f(x^0)$. 

\begin{lemma} \label{lem:Lnv}
Let $L_n \in \R^{n \times n}$ be defined by
\begin{align*}
    [L_n]_{i,i}&=\frac{12(3n-2)}{\de_i^{2}(3n+1)}, \quad i \in \{1, 2, \dots, n\}, \\
    [L_n]_{i,j}&= \frac{-36}{\de_i \de_j(3n+1)}, \quad i,j \in \{1, 2, \dots, n\}, i \neq j.
\end{align*}
Let $v \in \R^n$ be defined by
    $$
    v_i=\nintrecc x_i \nabla f(x^0)^\top x\dexxe, \quad i \in \{1, 2, \dots, n\}.
    $$
Then
 $\frac{1}{\de} L_n v=\nabla f(x^0).$
\end{lemma}
\begin{proof}
 First, we find an expression for $v_i, i \in \{1, 2, \dots, n\}.$ To make notation tighter, let $g=\nabla f(x^0).$ We have
\begin{align*}
    v_i&= \nintrecc x_i g^\top x \dexxe \\
    &=\nintrecc x_i \left ( \sum_{j=1}^n g_j x_j \right ) \dexxe \\
    &=\nintrecc x_i^2 g_i \dexxe + \sum_{j \neq i} \nintrecc x_i x_j g_j \dexxe \\
    &=\frac{\Delta_i^3}{3}\frac{\Delta}{\Delta_i}g_i+\sum_{j \neq i} \frac{\Delta_i^2}{2} \frac{\Delta_j^2}{2} \frac{\Delta}{\Delta_i \Delta_j} g_j \\
    &=\frac{\de_i^{2} \de}{3} g_i + \sum_{j \neq i} \frac{\de_i \de_j \de}{4} g_j=\frac{\de \de_i}{12} \left ( 4\de_i g_i+3 \sum_{j \neq i}  \de_j g_j \right ).
\end{align*}
Let $s=\sum_{j=1}^n \de_j g_j.$ Then
\begin{align*}
    v_i&=\frac{\de \de_i}{12} \left ( 4 \de_i g_i+3 \left (\sum_{j\neq i}^n \de_j g_j \right ) +3 \de_i g_i-3 \de_i g_i  \right )=\frac{\de \de_i}{12} \left (\de_i g_i+3s \right ).
\end{align*}
Now, let us compute $\frac{1}{\de} L_n v.$ Let $D \in \R^{n \times n}$ be the diagonal matrix with entries $\de_1, \de_2, \dots, \de_n.$ Note that 
\begin{align*}
    L_n&=\frac{12}{3n+1} D^{-1}L_n^\prime D^{-1},
\end{align*}
where $[L_n^\prime]_{i,i}=3n-2$ for all $i \in \{1, 2, \dots, n\}$ and $[L_n^\prime]_{i,j}=-3$ for all $i,j \in \{1, \dots, n\}, i \neq j.$ Let the vector $d= \bbm \de_1& \de_2&\cdots&\de_n \ebm^\top$ $\in \R^n.$ The vector $v$ can be written as
\begin{align*}
    v&=\frac{\de}{12}\left ( D^2g+3sd\right ).
\end{align*}
We obtain
\begin{align}
    \frac{1}{\de} L_n v&= \frac{12}{\de (3n+1)}D^{-1} L_n^\prime D^{-1}\frac{\de}{12} (D^2g+3sd) \notag \\
    &= \frac{1}{3n+1}D^{-1}  L_n^\prime D^{-1} (D^2g+3sd) \notag\\
    &=\frac{1}{3n+1}D^{-1} L_n^\prime D g+ \frac{3s}{3n+1}D^{-1} L_n^\prime D^{-1}d. \label{eq:Lnv}
\end{align}
The first term in \eqref{eq:Lnv} is equal to
\begin{align*}
    \frac{1}{3n+1}D^{-1} L_n^\prime D g
    &=\frac{1}{3n+1} D^{-1}\bbm (3n-2)\de_1 g_1 -3 \sum_{j \neq 1} \de_j g_j \\ \vdots \\ (3n-2) \de_n g_n -3 \sum_{j \neq n}\de_j g_j \ebm  \\
    &=\frac{1}{3n+1} D^{-1}\bbm (3n-2)\de_1 g_1 -3 \sum_{j \neq 1} \de_j g_j -3 \de_1 g_1 +3 \de_1 g_1 \\ \vdots \\ (3n-2) \de_n g_n -3 \sum_{j \neq n}\de_j g_j -3 \de_n g_n + 3 \de_n g_n \ebm  \\
    &=\frac{1}{3n+1} D^{-1}\bbm (3n+1)\de_1 g_1 -3s \\ \vdots \\(3n+1)\de_n g_n -3s\ebm  \\
    &= g-\frac{3s}{3n+1}D^{-1} \ones_n. 
\end{align*}
The second term in \eqref{eq:Lnv} is equal to 
\begin{align*}
    \frac{3s}{3n+1}D^{-1}L'_n D^{-1} d&=\frac{3s}{3n+1}D^{-1}L'_n \ones_n  \\
    &=\frac{3s}{3n+1} D^{-1} \bbm (3n-2)-3(n-1) \\ \vdots \\(3n-2)-3(n-1) \ebm  \\
    &= \frac{3s}{3n+1} D^{-1} \ones_n . 
\end{align*}
Hence, 
\begin{equation*}
    \frac{1}{\de} L_n = \frac{1}{3n+1}D^{-1} L'_n D g+ \frac{3s}{3n+1}D^{-1} L'_n D^{-1}d =g-\frac{3s}{3n+1} D^{-1} \ones_n+\frac{3s}{3n+1} D^{-1} \ones_n=g.\qedhere
\end{equation*}
\end{proof}

By the previous lemma, we now know that the error bound ad infinitum of the GSG is
\begin{align}\label{eq:ebsplit}
    \left  \Vert \limnvec \nabla_s f(x^0;S) -\nabla f(x^0) \right  \Vert&= \left \Vert \frac{1}{\de}L_n v+\frac{1}{\de}L_n w -\nabla f(x^0)\right \Vert= \frac{1}{\de} \left \Vert  L_n w \right \Vert.
\end{align}
To determine an upper bound for $\frac{1}{\de} \Vert L_n w \Vert$, recall that $L_n$ can be written as $D^{-1} \ddot{L}_n D^{-1},$ where  the diagonal matrix $D=\diag \bbm \de_1&\de_2&\cdots&\de_n\ebm \in \R^{n \times n}$ and the entries of $\ddot{L}_n \in \R^{n \times n}$ are given by
\begin{align}
    [\ddot{L}_n]_{i,i}&=\frac{12(3n-2)}{3n+1}, \quad i \in \{1, \dots, n\}, \notag \\
    [\ddot{L}_n]_{i,j}&=\frac{-36}{3n+1}, \quad i,j \in \{1, 2, \dots, n\}, i \neq j. \label{eq:ddotln}
\end{align}
Hence, the right-hand side of \eqref{eq:ebsplit} is
\begin{align} \label{eq:normsplit}
    \frac{1}{\de} \left \Vert L_n w\right \Vert&=\frac{1}{\de} \Vert D^{-1} \ddot{L_n} D^{-1}w  \Vert \leq \frac{1}{\de} \Vert D^{-1} \Vert \Vert \ddot{L}_n \Vert \Vert D^{-1} w \Vert.
\end{align}
The norm of $D^{-1}$ is simply $1/\de_{\tt min}$ where $\de_{\tt min}=\min_{i \in \{1, \dots, n\}} \de_i$. In the following two lemmas, we find the value of the two other norms that appear in \eqref{eq:normsplit}, $\Vert \ddot{L}_n \Vert$ and $\Vert D^{-1} w \Vert.$ 

\begin{lemma} \label{lem:normddotLn}
Let  $\ddot{L}_n \in \R^{n \times n}$ be defined as in \eqref{eq:ddotln}.
Then
\begin{align*}
 \left \Vert  \ddot{L}_n \right \Vert &=12. 
\end{align*}
\end{lemma}

\begin{proof}
 We find the norm by finding the largest singular value of $\ddot{L}_n^\top\ddot{L}_n.$ We have
\begin{equation*} \label{eq:normLn}
   \left [\ddot{L}_n^\top \ddot{L}_n \right ]_{i,j}=\left\{
  \begin{aligned}
  &\frac{144}{(3n+1)^2}(9n^2-3n-5), && \text{if}\, i=j, \\
    &\frac{144}{(3n+1)^2} (-9n-6), && \text{if}\, i \neq j, i,j \in \{1, 2, \dots, n\}.
  \end{aligned}\right.
\end{equation*} 
Let $t=\frac{144}{(3n+1)^2}.$
It follows that
\begin{align} \label{eq:LntLn}
\vert \ddot{L}_n^\top \ddot{L}_n -\lambda \id \vert &=t^n \left \vert \bbm 9n^2-3n-5-\frac{\lambda}{t}&-9n-6&\cdots&-9n-6&-9n-6 \\ -9n-6&9n^2-3n-5-\frac{\lambda}{t}&\cdots&-9n-6&-9n-6\\ \vdots&\vdots&\ddots&\vdots&\vdots\\-9n-6&-9n-6&\cdots&9n^2-3n-5-\frac{\lambda}{t}&-9n-6 \\-9n-6&-9n-6&\cdots&-9n-6&9n^2-3n-5-\frac{\lambda}{t} \ebm\right \vert.
\end{align}
Now, we apply elementary row  and column operations on the matrix in \eqref{eq:LntLn} to make it an upper-diagonal matrix. First, let Row $i$ $=$ Row $i~-$ Row 1 for $i \in \{2, 3, \dots, n\}.$ Second, using the new matrix, let Column $1$ $=$ Column $1$ + $\sum_{i=2}^n$ Column $i$. This generates the matrix
\begin{align*}
     \bbm 1-\frac{\lambda}{t}&-9n-6&-9n-6&\cdots&-9n-6&-9n-6\\0&9n^2+6n+1-\frac{\lambda}{t}&0&\cdots&0&0\\0&0&9n^2+6n+1-\frac{\lambda}{t}&\cdots&0&0\\ \vdots&\vdots&\vdots&\ddots&\vdots&\vdots\\0&0&0&\cdots&9n^2+6n+1-\frac{\lambda}{t}&0\\0&0&0&\cdots&0&9n^2+6n+1-\frac{\lambda}{t}\ebm 
\end{align*}
and therefore we have 
\begin{align*}
    \vert\ddot{L}_n^\top \ddot{L}_n -\lambda \id \vert 
    &=t^n\left (1-\frac{\lambda}{t} \right )\left (9n^2+6n+1-\frac{\lambda}{t} \right )^{n-1}.
\end{align*}
Hence, the eigenvalues of $\ddot{L}_n^\top \ddot{L}_n$ are
\begin{align*}
    \lambda_1&=t = \frac{144}{(3n+1)^2} \quad \text{and} \quad \lambda_{2, 3, \dots, n}=\frac{t}{(9n^2+6n+1)} = 144.\\
\end{align*}
We see that the maximum eigenvalue, denoted $\lambda_{\tt max}(\ddot{L}_n^\top \ddot{L}_n),$ is $144$. 
Therefore,
\begin{align*}
    \left \Vert \ddot{L}_n \right \Vert&= \sqrt{\lambda_{\tt max}(\ddot{L}_n^\top \ddot{L}_n)}=12.
\end{align*}\end{proof}

\begin{lemma} \label{lem:normDw}
  Let $f: \dom f \subseteq \R^n \to \R$ be $\mathcal{C}^2$ on an open domain containing $R(x^0;d).$  Let $L_{\nabla f}$ denote the Lipschitz constant of $\nabla f$ on  $R(x^0;d).$ Let $D=\diag \bbm \de_1&\cdots&\de_n \ebm \in \R^{n \times n}$ and let $w \in \R^n$  be defined by
 \begin{align*}
     w_i&= \nintrecc x_i R_1 (x^0+x) \dexxe, \quad i \in \{1, 2, \dots, n\}, 
 \end{align*}
 where $R_1(x^0+x)$ is the remainder term of the first-order Taylor expansion of $f(x^0+x)$ about  $x^0$. Then
 \begin{align*}
    \Vert D^{-1}w \Vert&\leq  \frac{\sqrt{n}}{8} L_{\nabla f} \de \de_{S}^{2}.
 \end{align*}
 Moreover, if all $\de_i$ are equal (i.e.\ the sample region is a hypercube), then
 \begin{align*}
    \Vert D^{-1}w \Vert&\leq  \frac{\sqrt{n}}{24}\frac{(2n+1)}{n} L_{\nabla f} \de \de_{S}^{2}.
 \end{align*}
 \end{lemma}
 \begin{proof}
 First, let us find an upper bound for  $\left \vert \frac{w_i}{\de_i} \right \vert, i \in \{1, 2, \dots, n\}.$ We have
 \begin{align}
     \left \vert \frac{w_i}{\de_i} \right  \vert &= \left \vert \nintrecc \frac{x_i}{\de_i} R_1(x^0+x) \dexxe\right \vert \notag \\
     &= \frac{1}{2} \left \vert \nintrecc  \frac{x_i}{\de_i} x^\top \nabla^2 f(\xi) x \dexxe \right \vert  \notag \\
     &\leq \frac{1}{2} \nintrecc \left \vert \frac{x_i}{\de_i} \right \vert \Vert x \Vert^2 \Vert \nabla^2 f(\xi) \Vert \dexxe. \notag
 \end{align}
Using $\Vert \nabla^2 f(\xi) \Vert\leq L_{\nabla f}$ and $R(0;d) \subseteq \R^n_+$, \color{black} we obtain 
 \begin{align}
     &\frac{1}{2} \nintrecc \left \vert \frac{x_i}{\de_i} \right \vert \Vert x \Vert^2 \Vert \nabla^2 f(\xi) \Vert \dexxe \notag \\
     \leq&\frac{L_{\nabla f}}{2} \nintrecc \left ( \frac{x_i}{\de_i} \right ) \left ( \sum_{j=1}^n x_j^2 \right ) \dexxe \notag \\
     =&\frac{L_{\nabla f}}{2} \left ( \nintrecc \frac{x_i^3}{\de_i} \dexxe +  \sum_{j \neq i} \nintrecc \frac{x_i}{\de_i} x_j^2 \dexxe \right ) \notag\\
     =&\frac{L_{\nabla f}}{2} \left ( \frac{\de \de_i^2}{4}+ \sum_{j \neq i} \frac{\de \de_j^2}{6}\right ) \notag \\
     =& \frac{L_{\nabla f}}{24} \de \left ( 3\de_i^2+ 2 \sum_{j \neq i} \de_j^2 \right ) =\frac{L_{\nabla f}}{24} \de  \left ( \de_i^2+ 2\de_S^2\right ). \label{eq:wiineq}
 \end{align}
 Therefore,
 \begin{equation*} \label{eq:normw2}
     \Vert D^{-1}w \Vert \leq \sqrt{\sum_{i=1}^n \left ( \frac{L_{\nabla f}}{24} \de  \left ( \de_i^2+ 2\de_S^2\right )\right )^2 } 
     \leq \frac{L_{\nabla f}}{24}\de \sqrt{\sum_{i=1}^{n} (3\de_S^2)^2}=\frac{\sqrt{n}}{8}L_{\nabla f} \de \de_S^2.
 \end{equation*}
 When all $\de_i$ are equal, then \eqref{eq:wiineq}  becomes
 \begin{align*}
     \frac{L_{\nabla f}}{24} \de  \left ( \de_i^2+ 2\de_S^2\right )&= \frac{L_{\nabla f}}{24} \de  \left ( \frac{\de_S^2}{n}+ 2\de_S^2\right )
 \end{align*}
 and it follows that
 \begin{align*}
     \Vert D^{-1}w \Vert &\leq \frac{\sqrt{n}}{24} \frac{(2n+1)}{n} L_{\nabla f}\de \de_S^2. \qedhere
 \end{align*}
 \end{proof}
 We are now ready to introduce an error bound ad infinitum for the GSG.
\begin{theorem}[Error bound ad infinitum for the GSG] \label{thm:ebadinf}
 Let $f: \dom f \subseteq \R^n \to \R$ be $\mathcal{C}^2$ on an open domain containing $R(x^0;d)$ where $d=\bbm \de_1&\de_2&\cdots&\de_n \ebm>0$ and $x^0$ is the reference point. Let $\de_{S}$ be the radius of $S \in \R^{n \times N}$ as defined in \eqref{df:radius}. Let $\de_{\tt min}=\min_{i \in \{1, \dots, n\}} \de_i.$ Let $L_{\nabla f}$ denote the Lipschitz constant of $\nabla f$ on  $R(x^0;d).$ 
Then
\begin{align}
   \left  \Vert \limnvec \nabla_s f(x^0;S) -\nabla f(x^0) \right  \Vert &\leq  \frac{3}{2}\sqrt{n} \, L_{\nabla f}  \, \frac{\de_{S}^2}{\de_{\tt min}}. \label{eq:ebhyper} 
\end{align}
Moreover, if all $\de_i$ are equal, then 
\begin{align*}
    \left  \Vert \limnvec \nabla_s f(x^0;S) -\nabla f(x^0) \right  \Vert &\leq  \frac{1}{2}(2n+1) \, L_{\nabla f} \,\de_{S}. 
\end{align*}
\end{theorem}
\begin{proof} We have 
\begin{align*}
    \left \Vert \limnvec \nabla_s f(x^0;S)- \nabla f(x^0) \right \Vert &= \left \Vert \de^{-1} L_n(v+w)-\nabla f(x^0) \right \Vert \\
    &= \left \Vert \de^{-1} L_n v- \nabla f(x^0) + \de^{-1}L_n w \right \Vert \\
    &\leq \left \Vert \de^{-1} L_n v-\nabla f(x^0) \right \Vert + \de^{-1} \left \Vert D^{-1} \right \Vert \Vert \ddot{L}_n\Vert \Vert D^{-1}w \Vert.
\end{align*}
By Lemma \ref{lem:Lnv}, we know $\Vert \de^{-1} L_n v -\nabla f(x^0) \Vert=0.$ By Lemma \ref{lem:normddotLn}, Lemma \ref{lem:normDw} and since $\Vert D^{-1} \Vert=\de_{\tt min},$ we obtain
\begin{align*}
   \left  \Vert \limnvec \nabla f(x^0;S) -\nabla f(x^0) \right  \Vert&= \de^{-1} \left \Vert D^{-1} \right \Vert \Vert \ddot{L}_n\Vert \Vert D^{-1}w \Vert \leq \frac{3}{2}\sqrt{n} \, L_{\nabla f} \, \eta \, \de_{S}.
\end{align*}
When all $\de_i$ are equal, $\de_{\tt min}=\de_i=\de_S/\sqrt{n}$ for any $i \in \{1, \dots, n\}.$ We obtain
\begin{align*}
    \left  \Vert \limnvec \nabla f(x^0;S) -\nabla f(x^0) \right  \Vert &\leq \de^{-1} \left \Vert D^{-1} \right \Vert \Vert \ddot{L}_n\Vert \Vert D^{-1}w \Vert \\
    &\leq \de^{-1} \frac{\sqrt{n}}{\de_S} (12)  \frac{\sqrt{n}}{24} \frac{(2n+1)}{n} L_{\nabla f}\de \de_S^2 \\
    &= \frac{1}{2}(2n+1) \, L_{\nabla f} \,\de_{S}.
\end{align*}
\end{proof}
In the previous theorem, we see that the error bound is $O(\frac{\de_{S}^2}{\de_{\tt min}})$.  As $\de_{S}$ is the radius of the sample region and $\de_{\tt min}$ is the length of the shortest side of the sample region, the theorem suggests that the more uniform the sample region the smaller the error.  In other words, we want the simplex with vertices $x^0, x^0+\Delta_1, \dots, x^0 +\Delta_n$ to be `as uniform as possible'. Analyzing the ratio $\Delta_S/\Delta_{\tt min},$ we see that the minimum value of this ratio is $\sqrt{n}$, which occurs when the hyperrectangle is a hypercube. 

\section{The  GSG over a  ball} \label{sec:ball}
In this section, we find an error bound ad infinitum for the GSG of $f$ at $x^0$ over  a ball. First, we present some results on integration over a ball. In the next theorem,  $P(x)$ denotes a monomial. That is,
\begin{align}\label{eq:monomial}
    P(x)&=x^\alpha=x_1^{\alpha_1} x_2^{\alpha_2} \cdots x_n^{\alpha_n}, 
\end{align}
where $\alpha_i \in \N\cup\{0\}$ for all $i.$ 

\begin{theorem}[Integrating over a ball]\emph{\cite{folland2001}} \label{thm:intball}
Let $P$ be a monomial defined as in \eqref{eq:monomial}. Let $\sigma$ be the $(n-1)$-dimensional surface measure on $S_n(0;r).$ Let $\beta_i=\frac{1}{2}(\alpha_i+1)$ for all $i$. Then
\begin{align*}
    \intballzero P(x)dx&=\frac{r^{\alpha_1+\dots+\alpha_n+n}}{\alpha_1+\dots+\alpha_n+n} \int_{S_n(0;r)}P d\sigma,
\end{align*}
where
\begin{equation*}
  \int_{S_n(0;r)} P d\sigma = \left\{
  \begin{aligned}
  &0, && \text{if any $\alpha_i$ is odd,} \\
    & \frac{2 \Gamma({\beta_1}) \Gamma(\beta_2)\cdots \Gamma(\beta_n)}{\Gamma(\beta_1+\beta_2+\dots+\beta_n)}, && \text{if all $\alpha_i$ are even}. \\
  \end{aligned}\right.
\end{equation*} 
\end{theorem}
We will also need an expression for the integral of $Q(x)=\vert x_1 \vert^{\alpha_1} \vert x_2 \vert^{\alpha_2} \cdots \vert x_n \vert^{\alpha_n}$ over the  ball $B_n(0;r).$
\begin{proposition} \label{prop:integrateabsvalue}
Let $Q(x)=\vert x_1 \vert^{\alpha_1} \vert x_2 \vert^{\alpha_2} \cdots \vert x_n \vert^{\alpha_n}$, and let $\beta_i=\frac{1}{2}(\alpha_i+1)$ for all $i$. Then 
\begin{align*}
    \intballzero  Q(x)dx&=\frac{2r^{\alpha_1+\dots+\alpha_n+n}}{(\alpha_1+\dots+\alpha_n+n)} \frac{\Gamma(\beta_1)\cdots\Gamma(\beta_n)}{\Gamma(\beta_1+\dots+\beta_n)}.
\end{align*}
\end{proposition}
\begin{proof}
Note that $\vert x_i \vert^{\alpha_i}$ is an even function for any $\alpha_i \in \N\cup\{0\}.$ Using this fact and following the same scheme of the proof for Theorem \ref{thm:intball} in \cite{folland2001} yields the result.
\end{proof}

Now, let us define the matrix of directions $S_R$ that is used to form the sample points. Recall that in $\R^n$, the conversion from Cartesian coordinates $x= \bbm x_1&x_2&\ldots&x_n \ebm^\top$ to $n$-spherical coordinates is
\begin{align*}
x_1&=\rho\cos\phi_1,\\
x_2&=\rho\sin\phi_1\cos\phi_2,\\
x_3&=\rho\sin\phi_1\sin\phi_2\cos\phi_3,\\
&\qquad\vdots\\
x_{n-2}&=\rho\sin\phi_1\cdots\sin\phi_{n-3}\cos\phi_{n-2},\\
x_{n-1}&=\rho\sin\phi_1\cdots\sin\phi_{n-2}\cos\theta,\\
x_n&=\rho\sin\phi_1\cdots\sin\phi_{n-2}\sin\theta,
\end{align*}where $\rho=\|x\|$ and $\theta, \phi_1, \dots, \phi_{n-2}$ are the angles that identify the direction of $x$. The angles have domains $\theta\in[0,2\pi)$ and $\phi_i\in[0,\pi)$ for all $i\in\{1,2,\ldots,n-2\}$. To keep the same notation as the hyperrectangle, we define
$$\Delta_1=r, \quad \Delta_2=2\pi, \quad \Delta_3=\Delta_4=\dots=\Delta_n=\pi.$$ As before, $N_i$ represents the number of subdivisions used to build the partitions in the ball. 
Once again, we define  $\overline{\Delta}_i=\Delta_i/N_i$ for all $i\in \{1, 2, \dots, n\},$  $\overline{\Delta}=\overline{\Delta}_1 \overline{\Delta}_2 \cdots \overline{\Delta}_n,$ and $N=N_1 N_2 \cdots N_n.$

Now we build the matrix $S_R \in \R^{n \times N}.$ The matrix $S_R$ contains all directions to add to the reference point $x^0$ to obtain a sample point in each partition of the ball $B_n(x^0;r)$. A polar grid is built, in which each partition is a ``bent'' hyperrectangle.  When using $S_R,$ the sample point chosen in each partition is the  rightmost  endpoint (the point with the largest values of $\rho, \theta, \phi_1, \dots, \phi_{n-2}$). Let $\vv{y}=\bbm y_1&y_2&\ldots&y_n\ebm^\top$ be a vector of indices in $\R^n$ (not $\R^{n-1}$ as it is the case for $\vv{z}$). Define 
\begin{align*}
    s^{\vv{y}}&=\frac{\rho y_1}{N_1} \bbm \cos{\frac {\pi y_3}{N_3}}& \sin{\frac{\pi y_3}{N_3}}\cos{\frac{\pi y_4}{N_4}}& \dots& \sin{\frac{\pi y_3}{N_{3}}} \cdots \sin{\frac{\pi y_n}{N_{n}}} \cos{\frac{2\pi y_2}{N_2}}& \sin{\frac{\pi y_3}{N_3}} \cdots \sin{\frac{\pi y_n}{N_n}} \sin{\frac{2\pi y_2}{N_2}}\ebm^\top \in \R^n.
\end{align*}
Then $S_R$ can be written as
\begin{align*}
    S_R&=\bbm s^{1, 1, 1, \dots,1, 1}&s^{1, 1, 1,  \dots, 1, 2}&\cdots&s^{N_1, N_2, N_3, \cdots, N_{n-1}, N_n} \ebm.
\end{align*}

\begin{figure}[H]
\centering
\includegraphics[width=0.65\textwidth]{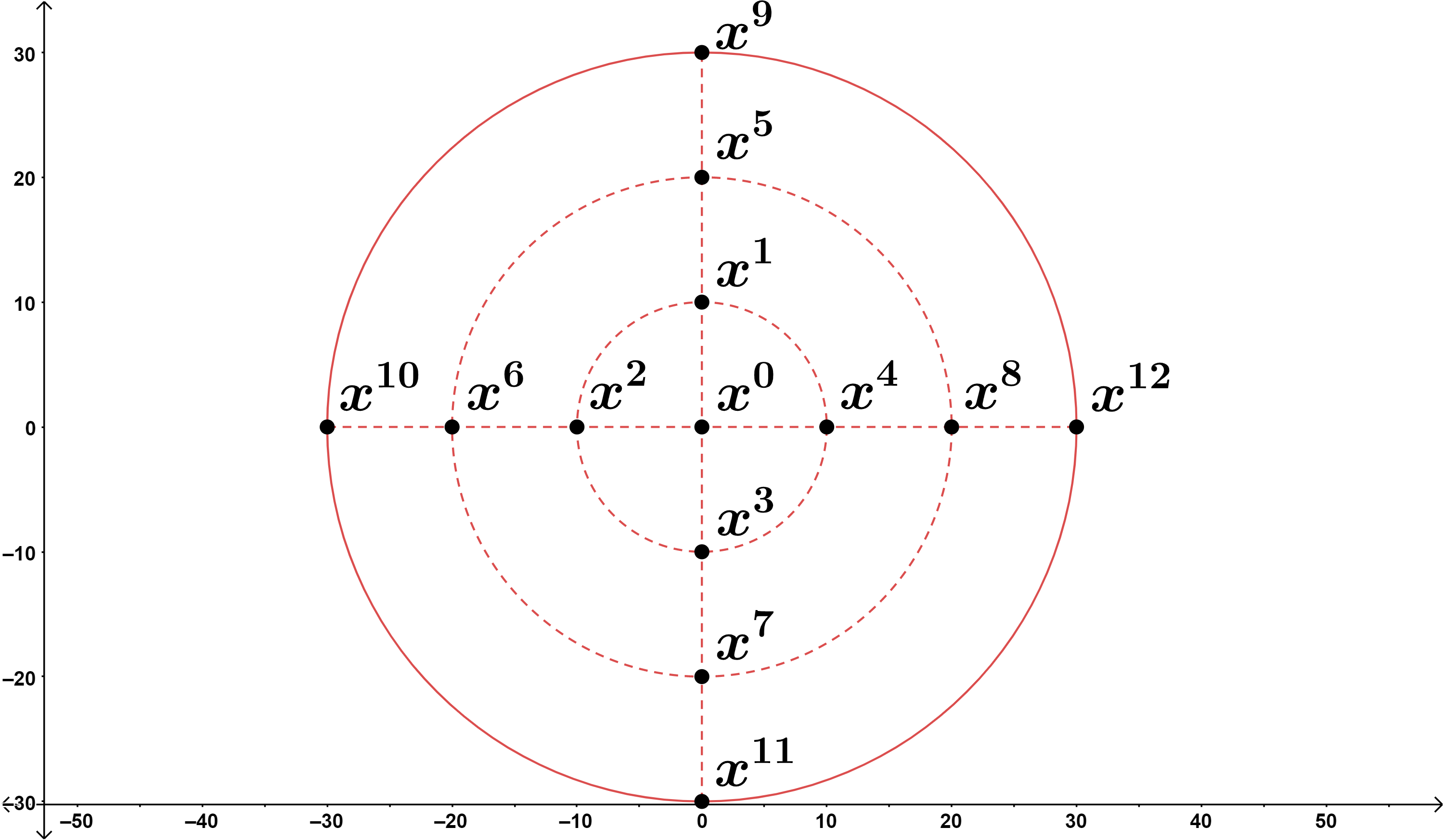}
\caption{An example of a sample set built from $S_R$ in $\R^2$.}
 \label{fig:SRball}
\end{figure}

Let us provide an example of a sample set built by using $S_R$ in $\R^2.$
\begin{example}
In this example, the reference point is $x^0=\bbm 0&0\ebm^\top.$ The sample region is a ball with radius $\Delta_1=r=30.$ Set $N_1=3$ and $N_2=4.$ Hence, $\overline{\Delta}_1=\frac{30}{3}=10$ and $\overline{\Delta}_2=\frac{2\pi}{4}=\frac{\pi}{2}.$ The matrix $S_R \in \R^{2 \times 12}$ is given by
\begin{align*}
    S_R&=\bbm s^{1,1}& s^{1,2}& s^{1,3}& s^{1,4}& s^{2,1}& s^{2,2}& s^{2,3}& s^{2,4} & s^{3,1} & s^{3,2}&s^{3,3}&s^{3,4}\ebm \\
    S_R&=\bbm 0&-10&0&10&0&-20&0&20&0&-30&0&30\\10&0&-10&0&20&0&-20&0&30&0&-30&0\ebm
\end{align*}
The sample points $x^j$ are obtained by  setting  $x^j=x^0+S_R e_j$ for all $i \in \{1, 2, \dots, 12\}.$ Figure \ref{fig:SRball} illustrates the sample points built from the matrix $S_R.$ 
\end{example}
Note that $S_R$ is full row rank whenever all $N_i>2$. For the remainder of this section, assume $N_i>2$ for all $i$. Hence, we want to find the limit  of the following expression:
\begin{align}
    \limnvec \nabla_s f(x^0;S_R)&=\limnvec \left (S_RS_R^\top \right )^{-\top} S_R \delta_f (x^0;S_R). \label{eq:limittosolve}
\end{align}

\noindent Define the determinant of the Jacobian as a function $J:\R^{n} \to \R$: 
\begin{align}
    J(y_1,\dots, y_n)&=\left (\frac{\rho y_1}{N_1} \right )^{n-1} \sin^{n-2}\frac{\pi y_3}{N_3}\sin^{n-3}\frac{\pi y_4}{N_4}\cdots\sin^2\frac{\pi y_{n-1}}{N_{n-1}}\sin \frac{\pi y_n}{N_{n}}. \label{eq:detjacobian}
\end{align}
  Let $J \in \R^{N \times N}$ be the matrix defined by
\begin{align*}
    J&=\diag\bbm j^{1,1,\dots,1}&j^{1, 1, \dots, 2}&\dots& j^{N_1, N_2,\dots, N_n}\ebm
\end{align*}
where $j^{\vv{y}}=J(y_1, y_2, \dots,y_n).$

Define
\begin{align*}
    K&=S_R J S_R^\dagger \\
    &=S_R J S_R^\top (S_RS_R^\top)^{-1}.
\end{align*}
Notice that  $K \in \R^{n \times n}$ is an invertible diagonal matrix such that 
\begin{align}\label{eq:K}
    S_R J= K S_R.
\end{align}

\noindent Considering \eqref{eq:limittosolve}, notice that
\begin{align*}
    (S_RS_R^\top)^{-\top}S_R \delta_f(x^0;S_R)
    &=(S_RS_R^\top)^{-\top}\frac{1}{\overline{\Delta}}K^{-1} K S_R \delta_f(x^0;S_R) \overline{\Delta} \\
    &=(K S_R S_R^\top\overline{\Delta})^{-\top} K S_R \delta_f(x^0;S_R)\overline{\Delta}\\
    &=(S_RJS_R^\top \overline{\Delta})^{-\top}S_R J\delta_f(x^0;S_R)\overline{\Delta}.
\end{align*}
Therefore, \eqref{eq:limittosolve} is equal to
\begin{align*}
    &\limnvec \left (S_RS_R^\top \right )^{-\top} S_R \delta_f (x^0;S_R) \\
    =&\limnvec  \left [\left (\sum_{y_1=1}^{N_1} \cdots \sum_{y_n=1}^{N_n} s^{\vv{y}} (s^{\vv{y}})^\top \right )^{-\top} \sum_{y_1=1}^{N_1} \cdots \sum_{y_n=1}^{N_n} s^{\vv{y}} \delta_f(x^0;s^{\vv{y}}) \right ] \\
    =&\limnvec  \left [\left (\sum_{y_1=1}^{N_1} \cdots \sum_{y_n=1}^{N_n} s^{\vv{y}} (s^{\vv{y}})^\top J(y_1, \dots, y_n)\overline{\Delta} \right )^{-\top} \sum_{y_1=1}^{N_1} \cdots \sum_{y_n=1}^{N_n} s^{\vv{y}} \delta_f(x^0;s^{\vv{y}}) J(y_1,\dots, y_n) \overline{\Delta} \right ].
\end{align*}
Considering 
    \begin{align}\limnvec  \left (\sum_{y_1=1}^{N_1} \cdots \sum_{y_n=1}^{N_n} s^{\vv{y}} (s^{\vv{y}})^\top J(y_1, \dots, y_n)\overline{\Delta} \right )^{-\top}
    ~\mbox{and}~ \limnvec \left ( \sum_{y_1=1}^{N_1} \cdots \sum_{y_n=1}^{N_n} s^{\vv{y}} \delta_f(x^0;s^{\vv{y}}) J(y_1, \dots, y_n) \overline{\Delta} \right ),
    \label{eq:twolimitssphere}\end{align}
we shall show both limits exist, so can write
\begin{align*}
    &\limnvec \left (S_RS_R^\top \right )^{-\top} S_R \delta_f (x^0;S_R) \notag \\
    =&\limnvec  \left (\sum_{y_1=1}^{N_1} \cdots \sum_{y_n=1}^{N_n} s^{\vv{y}} (s^{\vv{y}})^\top J(y_1, \dots, y_n)\overline{\Delta} \right )^{-\top} \, \limnvec \left ( \sum_{y_1=1}^{N_1} \cdots \sum_{y_n=1}^{N_n} s^{\vv{y}} \delta_f(x^0;s^{\vv{y}}) J(y_1, \dots, y_n) \overline{\Delta} \right ).
\end{align*}

\noindent We begin by examining the first limit in \eqref{eq:twolimitssphere}.  Recall the following formula for the volume of a ball in $\R^{n+2}$:
\begin{align}\label{eq:volumenplus2}
    V_{n+2}&=\frac{2 \pi^{\frac{n}{2}+1}r^{n+2}}{\Gamma(\frac{n}{2}+1)(n+2)}, 
\end{align}
where $\Gamma$ is the Gamma function given by 
\begin{equation*} 
  \Gamma \left (\frac{n}{2}+1 \right ) \in \left\{
  \begin{aligned}
  &(n-1)! && \text{if}\quad \frac{n}{2} \in \N, \\
    &\left (\frac{n}{2}-1 \right )\left (\frac{n}{2}-2 \right )\cdots \frac{1}{2} \sqrt{\pi} && \text{if}\quad \frac{n}{2} \notin \N. \\
  \end{aligned}\right.
\end{equation*} 
 Let 
 \begin{align}\label{eq:matrixMn}
     M_n&=\limnvec \left ( \sum_{y_1=1}^{N_1} \cdots \sum_{y_n=1}^{N_n} s^{\vv{y}} \left (s^{\vv{y}} \right )^\top J(y_1, \dots, y_n) \overline{\Delta} \right ) \in \R^{n \times n}.
 \end{align}
 
\noindent Proposition \ref{prop:firstlimitball} finds the matrix $M_n$ and thereby provides the first limit in \eqref{eq:twolimitssphere}.

\begin{proposition}\label{prop:findingMn} Let $M_n \in \R^{n \times n}$ be defined as in \eqref{eq:matrixMn}. Then 
\begin{equation*} 
 [M_n]_{i,j}= \left\{
  \begin{aligned}
  & \frac{V_{n+2}}{2\pi}&& \text{if} \quad  i=j, \\
    &0 && \text{if}\quad i \neq j. \\
  \end{aligned}\right.
\end{equation*} 
Consequently \label{prop:firstlimitball}
\begin{align*}
    \limnvec \left (\sum_{y_1=1}^{N_1} \cdots  \sum_{y_n=1}^{N_n} s^{\vv{y}} (s^{\vv{y}})^\top J(y_1, \dots, y_n)\overline{\Delta} \right )^{-\top}&=\frac{2\pi}{V_{n+2}} \id_n.
\end{align*}
\end{proposition}
\begin{proof}
Each entry of $M_n$ is a  $n$-tuple Riemann sum  over $B_n(0;r).$ Taking the limit as $\vv{N} \to \vv{\infty},$ each entry of $M_n$ can be written as the following integral (in Cartesian coordinates):
 \begin{align*} 
     [M_n]_{i,j}&= \intballzero x_ix_j  dx,  \quad i,j \in \{1, 2, \dots, n\}. 
 \end{align*}
The off-diagonal entries of $M_n$ are zero by Theorem \ref{thm:intball}. The diagonal entries are given by
\begin{align*}
    [M_n]_{i,i}&=\frac{r^{n+2}}{(n+2)} \frac{2 \Gamma \left (\frac{1}{2} \right )^{n-1}\Gamma \left (\frac{3}{2} \right )}{\Gamma \left ( \frac{1}{2} (n-1)+\frac{3}{2} \right )} =\frac{r^{n+2}}{(n+2)}\frac{2\pi^{\frac{n-1}{2}} \pi^{\frac{1}{2}}}{2 \Gamma \left ( \frac{n}{2}+1 \right )}=\frac{r^{n+2} \pi^{\frac{n}{2}}}{(n+2) \Gamma \left ( \frac{n}{2}+1\right )}.
\end{align*}
From \eqref{eq:volumenplus2}, we see that the diagonal entries are simply
\begin{align*}
    [M_n]_{i,i}&=\frac{V_{n+2}}{2\pi}.
\end{align*}
The second equation follows trivially.\qedhere
\end{proof}

\noindent In the next proposition, we give an expression for the second limit in \eqref{eq:twolimitssphere}.
\begin{proposition}\label{prop:secondlimitball}
Let $f:\dom f \subseteq \R^n \to \R$ with $B_n(x^0;r) \subseteq \dom f.$ Then the following limit  can be written as a vector of integrals (in Cartesian coordinates):
\begin{align} \label{eq:Tnball}
    \limnvec \left [ \sum_{y_1=1}^{N_1} \cdots \sum_{y_n=1}^{N_n} s^{\vv{y}} \delta_f(x^0;s^{\vv{y}}) J(y_1, \dots, y_n) \overline{\Delta} \right ]&= \bbm \intballzero x_1 \left ( f(x^0+x)-f(x^0) \right ) dx\\  \intballzero  x_2 \left ( f(x^0+x)-f(x^0) \right )  dx \\ \vdots \\  \intballzero  x_n \left ( f(x^0+x)-f(x^0) \right ) dx\ebm =T_n \in \R^n.
\end{align}
\end{proposition}
\begin{proof}
The $n$-tuple sum of the left-hand side of \eqref{eq:Tnball} is a Riemann sum  with  a finite-sized  sample region  $B_n(x^0;r).$ The result follows by taking the limit as $\vv{N} \to \vv{\infty}.$
\end{proof}
\noindent Now we generalize the matrix $S_R$. Let
\begin{align*}
    S&=\bbm s^{1, 1, 1, \dots,1, 1}&s^{1, 1, 1,  \dots, 1, 2}&\cdots&s^{N_1, N_2, N_3, \cdots, N_{n-1}, N_n} \ebm
\end{align*}
 be a matrix in $\R^{n \times N}$ in which each column $s$ is a direction to add to $x^0$ to form exactly one arbitrary sample point in each partition of $B_n(x^0;r).$ 
Note that Propositions \ref{prop:firstlimitball} and \ref{prop:secondlimitball} still hold by considering $S$ instead of $S_R.$ Indeed, since $f$ is a continuous function, as long as exactly  one sample point is used in each partition of the ball, the results of the previous two propositions hold. We are now ready to provide an expression for the GSG ad infinitum of $f$ at $x^0$ over $B_n(x^0;r).$ 

\begin{theorem}[The GSG over a ball] \label{thm:limitgsgball}
Let $f:\dom f \subseteq \R^n \to \R$ with $B_n(x^0;r) \subseteq \dom f.$ Let $S \in \R^{n \times N}$ be a matrix such that each sample point $x^0+Se_j, j \in \{1, 2, \dots, N\},$ is in exactly one partition of the ball $B_n(x^0;r)$. Let $V_{n+2}$ be the volume of a ball with radius $r$ in $\R^{n+2}$ and  $T_n \in \R^n$ be defined as in \eqref{eq:Tnball}. Then
\begin{align}
    \limnvec \nabla_s f(x^0;S)&= \frac{2\pi}{V_{n+2}}T_n. \label{eq:limgsgball}
\end{align}
\end{theorem}
\begin{proof}
we have
\begin{align*}
    &\limnvec \nabla_s f(x^0;S)\\
    &\quad=\limnvec \left (S S ^\top \right )^{-\top} S \delta_f (x^0;S) \\
    &\quad=\limnvec  \left [\left (\sum_{y_1=1}^{N_1} \cdots \sum_{y_n=1}^{N_n} s^{\vv{y}} (s^{\vv{y}})^\top J(y_1, \dots, y_n)\overline{\Delta} \right )^{-\top} \right ] \limnvec \left [ \sum_{y_1=1}^{N_1} \cdots \sum_{y_n=1}^{N_n} s^{\vv{y}} \delta_f(x^0;s^{\vv{y}}) J(y_1, \dots, y_n) \overline{\Delta} \right ] \\
    &\quad=\frac{2\pi}{V_{n+2}} T_n,
\end{align*}
by Propositions \ref{prop:firstlimitball} and \ref{prop:secondlimitball}.
\end{proof}
\noindent Let us provide an example of the calculations necessary to obtain the limit of the the GSG.
\begin{example} \label{ex:exampleball}
Let $f:\R^2 \to \R: x \mapsto x_1^2+x_2^2.$ Let the reference point be $x^0=\bbm3&1\ebm^\top.$ Let the sample region be $B_2(x^0;1).$ By Theorem \ref{thm:limitgsgball}, we know that
\begin{align*}
    \limnvec \nabla_s f(x^0;S)&=\frac{2\pi}{V_4} \bbm \int_{B_2(0;1)} x_1 \left (f(x^0+x)-f(x^0) \right ) dx \\ \int_{B_2(0;1)} x_2 \left (f(x^0+x)-f(x^0) \right ) dx\ebm.
\end{align*}
Writing the vector of integrals in polar coordinates and since $V_4=\frac{\pi^2}{4}$, we obtain
\begin{align*}
    \limnvec \nabla_s f(x^0;S)&=4\pi \bbm \int_{0}^{2\pi} \int_{0}^{1} r \cos{\theta} \left ( (3+r\cos{\theta})^2+(1+r \sin{\theta})^2-10 \right )r dr d\theta \\ \int_{0}^{2\pi} \int_{0}^{1} r \sin{\theta} \left ( (3+r\cos{\theta})^2+(1+r \sin{\theta})^2-10 \right )r dr d\theta\ebm=4\pi \bbm \frac{3}{2}\pi \\ \frac{\pi}{2} \ebm = \bbm 6 \\2 \ebm.
\end{align*}
Note that for this problem $\limnvec \nabla_s f(x^0;S)=\nabla f(x^0).$
\end{example}
The reason why the GSG is perfectly accurate in the previous example will be discussed at the end of this section, but first we develop an error bound ad infinitum for the GSG over a  ball. To obtain such an error bound, we require $f$  to be $\mathcal{C}^3$ on an open domain containing $B_n(x^0;r)$.  The function $f$ at  $x^0+x$ is written as the second-order Taylor expansion. That is,
 \begin{align}
     f(x^0+x)&=f(x^0)+\nabla f(x^0)^\top x+ \frac{1}{2} x^\top \nabla^2 f(x^0)x +R_2(x^0+x), \label{eq:secondordertaylor}
 \end{align}
 where the remainder term satisfies
 \begin{align*}
  \vert R_2(x^0+x) \vert &\leq  \frac{1}{6} L_H \Delta_S^3 
 \end{align*}
 and $L_H$ denotes the Lipschitz constant of the Hessian on $B_n(x^0;r)$.
By rewriting $f(x^0+x)$  in the form of \eqref{eq:secondordertaylor},  each entry of the vector $T_n$ defined in Proposition \ref{prop:secondlimitball} can be written as
\begin{align*}
    [T_n]_i&= \intballzero x_i \left (f(x^0+x)-f(x^0)\right ) dx \\
    &=\intballzero x_i \left ( \nabla f(x^0)^\top x+\frac{1}{2} x^\top \nabla^2 f(x^0) x +R_2(x^0+x) \right ) dx \\
    &=\intballzero x_i \nabla f(x^0)^\top x dx  + \intballzero x_i x^\top \nabla^2 f(x^0) x dx \intballzero x_i R_2(x^0+x) dx,
\end{align*}
for $i \in \{1, 2, \dots, n\}.$
Let $v \in \R^n$ be the vector defined by
\begin{align}
    v_i&=\intballzero x_i \nabla f(x^0)^\top x dx, \quad i \in \{1, 2, \dots, n\}, \label{eq:vectorvball}
\end{align}
$w \in \R^n$ be the vector defined by
\begin{align}
    w_i&=\intballzero x_i x^\top \nabla^2 f(x^0) x dx, \quad i \in \{1, 2, \dots, n\}, \label{eq:vectorwball}
\end{align}
and $z \in \R^n$ be the vector defined by 
\begin{align}
    z_i&=\intballzero x_i R_2(x^0+x) dx, \quad i \in \{1, 2, \dots, n\}. \label{eq:z}
\end{align}
Then the expression for $\limnvec \nabla_s f(x^0;S)$ given in \eqref{eq:limgsgball} can be written as
\begin{align}
    \limnvec \nabla_s f(x^0;S)&= \frac{2\pi}{V_{n+2}} (v+w+z)= \frac{2\pi}{V_{n+2}} v+ \frac{2\pi}{V_{n+2}} w + \frac{2\pi}{V_{n+2}}z.\label{eq:limgsg3terms}
\end{align}
We now find an expression for the first two terms in \eqref{eq:limgsg3terms}.

\begin{lemma} \label{lem:intv}
Let $f: \dom f \subseteq \R^n \to \R$ be $\mathcal{C}^3$ on an open domain containing $B_n(x^0;r).$  Let $v \in \R^n$ and $w \in \R^n$ be defined as in \eqref{eq:vectorvball} and \eqref{eq:vectorwball}.  Then
\begin{align*}
    \frac{2\pi}{V_{n+2}}v=\nabla f(x^0) ~\mbox{and}~    \frac{2\pi}{V_{n+2}}w=0. 
\end{align*}
\end{lemma}
\begin{proof}
Let $g_i=[\nabla f(x^0)]_i$ for all $i \in \{1, 2, \dots, n\}$. We have
\begin{align*}
    \frac{2\pi}{V_{n+2}} \intballzero x_i \nabla f(x^0)^\top x dx&=\frac{2\pi}{V_{n+2}} \intballzero x_i^2 g_i dx+\sum_{j \neq i }\intballzero x_ix_j g_i dx= \frac{2\pi}{V_{n+2}} \frac{V_n+2}{2 \pi} g_i+0=g_i,
\end{align*}
by Proposition \ref{prop:findingMn}. Therefore $\frac{2\pi}{V_{n+2}} v=\nabla f(x^0).$

Let $\nabla^2 f(x^0)=H \in \R^{n \times n}.$ We have 
\begin{align*}
    \frac{\pi}{V_{n+2}} \intballzero x_i x^\top \nabla^2 f(x^0) x  dx&=  \frac{\pi}{V_{n+2}} \left ( \sum_{j=1}^n \intballzero x_ix_j^2 H_{j,j} dx + 2 \sum_{j=1}^n \sum_{\substack{k=1\\ k \neq j}}^n \intballzero x_ix_jx_k H_{j,k} dx\right )=0,
\end{align*}
by Theorem \ref{thm:intball}. Therefore, $\frac{2\pi}{V_{n+2}}w=0.$
\end{proof}
\noindent Next, we find an upper bound for the third term in \eqref{eq:limgsg3terms}.  In Lemma \ref{lem:intR2} we create the redundant variable $\Delta_S = r$.  This allows for easy and immediate comparison to the results in Section \ref{sec:hyperrectangle}.
\begin{lemma} \label{lem:intR2}
 Let $f:\dom f \subseteq \R^n \to \R$ be $\mathcal{C}^3$ on an open domain containing $B_n(x^0;r).$ Let $z \in \R^n$ be defined as in \eqref{eq:z}. Let $$\eta=\frac{\Gamma(\frac{n+4}{2})}{\sqrt{\pi} \Gamma(\frac{n+3}{2})}.$$ Also, denote by $L_H$ the Lipschitz constant of the Hessian on $B_n(x^0;r).$  Let $\Delta_S = r$.
Then $$\frac{2\pi}{V_n+2} \Vert z \Vert \leq \frac{\sqrt{n}}{3 \sqrt{\pi}} L_H \eta \Delta_S^2.$$
\end{lemma}
\begin{proof}
We have
\begin{align}
    \vert z_i \vert &=\left \vert  \frac{2\pi}{V_{n+2}} \intballzero x_i R_2(x^0+x) dx \right \vert \notag \\
    &\leq \frac{2\pi}{V_{n+2}}\intballzero \vert x_i \vert \left \vert R_2(x^0+x) \right \vert dx \notag \\
    &\leq \frac{2\pi}{V_{n+2}} \frac{1}{6}  L_{H} \Delta_S^3 \intballzero \vert x_i \vert dx. \label{eq:intballabs}
\end{align}
By Proposition \ref{prop:integrateabsvalue}, we know that
\begin{align}
    \intballzero \vert x_i \vert dx =2 \frac{r^{n+1} \Gamma(1) \Gamma(\frac{1}{2})^{n-1}}{(n+1) \Gamma \left ( 1+\frac{n-1}{2}\right )} =\frac{2r^{n+1} \pi^{\frac{n-1}{2}}}{(n+1)\Gamma \left (\frac{n+1}{2}\right )}. \label{eq:intabsvalue}
\end{align}
Note that 
\begin{align*}
    \frac{2}{(n+1)\Gamma \left ( \frac{n+1}{2} \right )}&=\frac{1}{\Gamma \left ( \frac{n+3}{2}\right )}.
\end{align*}
Hence, \eqref{eq:intabsvalue} can be written as
\begin{align} \label{eq:answerintabs}
    \intballzero \vert x_i \vert dx &=\frac{\pi^{\frac{n-1}{2}}r^{n+1}}{\Gamma \left ( \frac{n+3}{2}\right )}=\frac{V_{n+1}}{\pi}.
\end{align}
Substituting \eqref{eq:answerintabs} in \eqref{eq:intballabs} gives
\begin{align}
    \left \vert  \frac{2\pi}{V_{n+2}} \intballzero x_i R_2(x^0+x) dx \right \vert &\leq \frac{1}{3} \frac{V_{n+1}}{V_{n+2}} L_H \Delta_S^3. \label{eq:ratiovolume}
\end{align}
The term  $V_{n+1}/V_{n+2}$ in \eqref{eq:ratiovolume} is
\begin{align*}
    \frac{V_{n+1}}{V_{n+2}}&=\frac{\pi^{\frac{n+1}{2}}r^{n+1}}{\Gamma \left (\frac{n+3}{2} \right )} \frac{\Gamma \left ( \frac{n+4}{2}\right )}{\pi^{\frac{n+2}{2}}r^{n+2}}=\frac{1}{\sqrt{\pi}\Delta_S} \eta.
\end{align*}
Thus,
\begin{align*}
   \vert z_i \vert  &\leq \frac{1}{3\sqrt{\pi}} \eta L_H \Delta_S^2, \quad \forall i \in \{1, 2, \dots, n\}.
\end{align*}
Therefore, \begin{align*}\Vert z \Vert  &\leq \frac{\sqrt{n}}{3 \sqrt{\pi}} L_H \eta \Delta_S^2.\qedhere\end{align*}
\end{proof}
\begin{theorem}[Error bound ad infinitum for the GSG over a ball] \label{thm:ebadinfball}
Let $f:\dom f \subseteq \R^n \to \R$ be $\mathcal{C}^3$ on an open domain containing $B_n(x^0;r).$ Let the vectors $v, w, z \in \R^n$ be defined as in \eqref{eq:vectorvball}, \eqref{eq:vectorwball}, \eqref{eq:z}, respectively. Let $$\eta=\frac{\Gamma(\frac{n+4}{2})}{\sqrt{\pi} \Gamma(\frac{n+3}{2})}.$$ Denote by $L_H$ the Lipschitz constant of $\nabla^2f$ on $B_n(x^0;r).$ Let   $\Delta_S$ be defined as in \eqref{df:radius}. Then
\begin{align} \label{eq:ebball}
    \left \Vert \limnvec \nabla_s f(x^0;S)-\nabla f(x^0) \right \Vert &\leq  \frac{\sqrt{n}}{3\sqrt{\pi}} L_{H} \eta \Delta_S^2.
\end{align}
\end{theorem}
\begin{proof}
We have
\begin{align*}
    \left \Vert \limnvec \nabla_s f(x^0;S)-\nabla f(x^0) \right \Vert&=\left \Vert  \frac{2\pi}{V_{n+2}}v+\frac{2\pi}{V_{n+2}}w+\frac{2\pi}{V_{n+2}}z -\nabla f(x^0)\right  \Vert \\
    &=\left \Vert \nabla f(x^0)+0+\frac{2\pi}{V_{n+2}}z-\nabla f(x^0) \right \Vert  \\
    &= \frac{2\pi}{V_{n+2}} \Vert z \Vert \leq \frac{\sqrt{n}}{3\sqrt{\pi}} L_{H} \eta \Delta_S^2.
\end{align*}
\end{proof}

Note that the error bound ad infinitum over a ball is $\mathcal{O}(\Delta_S^2)$, which is not the case for the error bound ad infinitum defined in Section \ref{sec:error}. This is due to the fact that, for each column $s \in S$, its opposite $-s$ is also in $S$ as $\vv{N} \to \vv{\infty}.$ Therefore, the limit of the GSG over $B_n(x^0;r)$ is equivalent to the limit of the generalized centered simplex gradient over a half-ball centered at $x^0$ of radius $r$. Notice that the shape of the sample region is not the key point to obtain an error bound  $\mathcal{O}(\Delta_S^2).$ The position of the reference point $x^0$ is what matters. Indeed, we could get an error bound ad infinitum  of accuracy $\mathcal{O}(\Delta_S^2)$ by considering a hyperrectangle, but instead of letting $x^0$ be the left endpoint of the sample region as we did in the previous sections, let $x^0$ be located at the intersection of all diagonals of the hyperrectangular sample region. 
Finally, note that the error bound in \eqref{eq:ebball} involves the Lipschitz constant of the Hessian of $f$, $L_H.$ Therefore, the error bound reduces to zero whenever $f$ is a polynomial of degree at most 2. This explains why the GSG is a perfect approximation in  Example \ref{ex:exampleball}.

\section{Conclusion and future research directions}\label{sec:conc}
In this paper, we have provided an expression for the GSG ad infinitum  and an error bound ad infinitum for the GSG over both a hyperrectangle and a ball.  In both cases, we note that the error bound is independent of the number of sample points, which is critical in allowing the analysis of the limits.  Examining the techniques used in each case, it seems likely that an error bound ad infinitum (independent of $N$) for the GSG of $f$ at $x^0$ over any reasonable sample region can be defined.  However, repeating the process for every possible region is clearly an unreasonable proposition.  A more practical open question is the following: Given a set of sample points $\Omega \subseteq \R^n$ and a bijection $\mathcal{T} : \R^n \mapsto \R^n$ such that $\mathcal{T}(\Omega) = R(x^0; d)$, can the bijection be used to determine the GSG ad infinitum  and an error bound ad infinitum for the GSG over $\Omega$?  

Comparing Theorems \ref{thm:ebadinf} and \ref{thm:ebadinfball}, we see that the position of the reference point $x^0$ has an impact on the error bound. Indeed, when $x^0$ is the center of the sample region, then we obtained an error bound ad infinitum of $\mathcal{O}(\Delta_S^2$).  But, when $x^0$ is on the boundary  of the sample region, then we obtained an error bound ad infinitum of $\mathcal{O}(\Delta_S)$.  It is unclear how these conclusions change if the reference point is at another location within the sample region.  

In \cite{MR4163088} it was found that under certain conditions, the limit of the classical error bound for the GSG in $\R^n$ (\eqref{eq:ebgsg} herein) could be taken directly. It is possible that this is true in $\R^n$ as well.  However, the techniques in \cite{MR4163088} do not adapt directly.

We conclude this paper with a comparison of classical error bounds (as given by 
\eqref{eq:ebgsg} and \eqref{eq:ebgcsg}) as $N$ gets large to the error bounds ad infinitum derived in Theorems \ref{thm:ebadinf} and \ref{thm:ebadinfball}. 
\begin{example}
In this example, we consider the function $f: \R^2 \to \R: x=\bbm x_1&x_2\ebm^\top\mapsto x_1^3+x_2^3.$ The reference point is set to $x^0=\bbm 1&1\ebm^\top.$ Two sample regions are considered: the square $[0,1]\times[0,1]$ and the ball $B_2(x^0;1).$ 

We set $N_1=N_2$, so $N=(N_1)^2$.  The classical error bounds are computed using \eqref{eq:ebgsg} and \eqref{eq:ebgcsg}.  The error bounds ad infinitum are computed using Theorems \ref{thm:ebadinf} and \ref{thm:ebadinfball}.  Finally, the GSG is constructed and the true absolute error is computed for both sample regions.  Figure \ref{fig:ebr2comparison} visualizes the results for $N_1 \in \{2^2, 2^3, \ldots, 2^{10}\}$ (so, $N \in \{2^4, 2^6, \ldots, 2^{20}\}$).   Note the  bounds ad infinitum are independent of $N_1$, so constants.

\begin{figure}[H]
\includegraphics[width=1\textwidth]{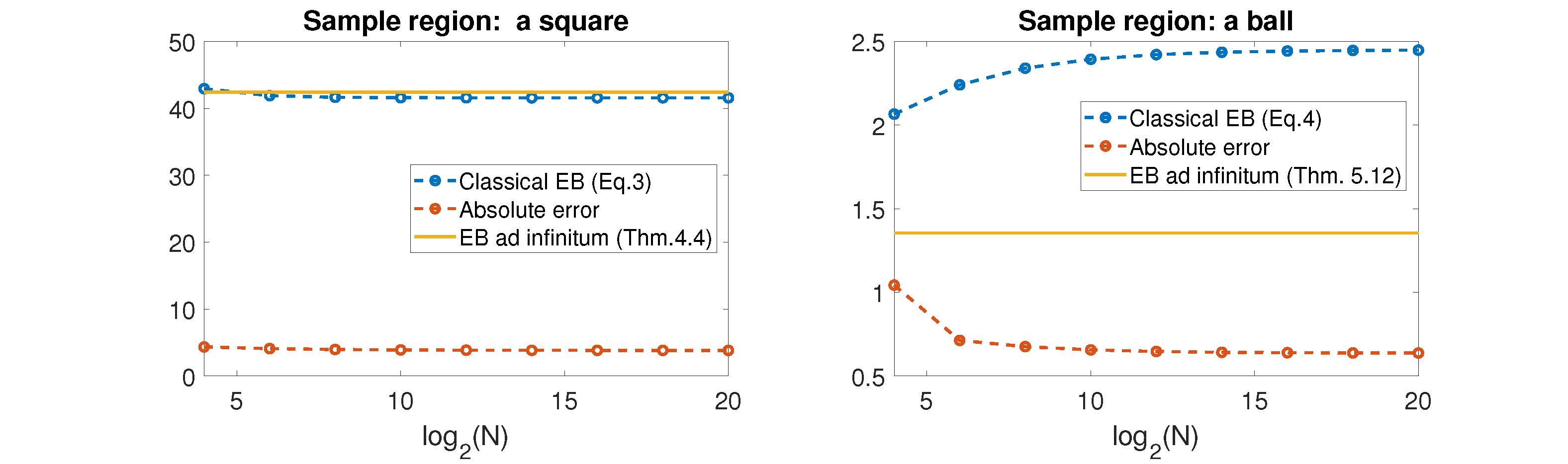} 
\caption{The error bound ad infinitum in $\R^2$ for two different sample regions}
\label{fig:ebr2comparison}
\end{figure}

Based on this example, the error bounds ad infinitum provides an accurate upper bound as low as $N=16$. It also appears that the error bound ad infinitum over the ball $B_2(x^0;1)$ provides a tighter error bound than the classical error bound for $N\geq 16$.  Finally, in this example, it appears that the classical error bounds may converge as $N$ tends to infinity as conjectured.
\end{example}

\section*{Acknowledgements}
Hare's research is partially funded by the Natural Sciences and Engineering Research Council (NSERC) of Canada, Discover Grant \#2018-03865.  Jarry-Bolduc's research is  partially funded by the Natural Sciences and Engineering Research Council (NSERC) of Canada, Discover Grant \#2018-03865. Jarry-Bolduc would like to acknowledge UBC for the funding received through the University Graduate Fellowship award.

\bibliographystyle{alpha}
\bibliography{DFOBIB}

\end{document}